\documentclass[12pt]{amsart}

\usepackage{stackrel}
\usepackage{latexsym}
\usepackage[pdftex]{graphicx}
\usepackage{setspace}

\newtheorem{theorem}{Theorem}
\newtheorem{lemma}{Lemma}
\newtheorem{proposition}{Proposition}
\newtheorem{remark}{Remark}

\newtheorem{corollary}{Corollary}

\def\R{{\mathbb R}}

\newcommand{\beq}{\begin{equation}}
\newcommand{\eeq}{\end{equation}}
\newcommand{\beqna}{\begin{eqnarray*}}
\newcommand{\eeqna}{\end{eqnarray*}}
\newcommand{\beqn}{\begin{equation*}}
\newcommand{\eeqn}{\end{equation*}}
\newcommand{\bp}{\begin{proof}}
\newcommand{\ep}{\end{proof}}
\newcommand{\bprop}{\begin{proposition}}
\newcommand{\eprop}{\end{proposition}}
\newcommand{\bt}{\begin{theorem}}
\newcommand{\et}{\end{theorem}}
\newcommand{\bex}{\begin{Example}}
\newcommand{\eex}{\end{Example}}
\newcommand{\bc}{\begin{corollary}}
\newcommand{\ec}{\end{corollary}}
\newcommand{\bl}{\begin{lemma}}
\newcommand{\el}{\end{lemma}}

\begin{document}

\title
[A problem of Klee]
{An asymmetric convex body with maximal sections of constant volume}

\author{Fedor Nazarov, Dmitry Ryabogin and Artem Zvavitch}
\address{Department of Mathematics, Kent State University,
Kent, OH 44242, USA} \email{nazarov@math.kent.edu}\email{ryabogin@math.kent.edu}\email{zvavitch@math.kent.edu}

\thanks{The first named author is supported in
part by U.S.~National Science Foundation Grant DMS-0800243. The second and third named authors are supported in part by U.S.~National Science Foundation Grant DMS-1101636. }



\begin{abstract}
We show that in all dimensions $d\ge 3$, there exists an asymmetric convex body of revolution all of whose maximal hyperplane sections
have the same volume. This gives the negative answer to the question posed by V. Klee in 1969.
\end{abstract}

\maketitle

\section{Introduction}

As usual,  a  {\it convex body} $K\subset {\mathbb R}^d$ is a compact convex subset  of ${\mathbb R}^d$ with non-empty interior.  We assume that the origin is contained in the interior of $ K$. We consider    the {\it maximal section function} $M_K$:
$$
M_K(u)=\max_{t\in {\mathbb R}}\mathrm{vol}_{d-1}(K\cap(u^{\perp}+tu)),\qquad u\in {\mathbb S}^{d-1}.
$$

 Here  $u^{\perp}$ stands for the hyperplane passing through the origin and orthogonal to the unit vector $u$, $K\cap (u^{\perp}+tu)$ is the section of $K$ by the affine hyperplane $u^{\perp}+tu$.
It is well known, \cite{Ga},  that for  origin-symmetric convex bodies the maximal sections are central (i.e., correspond to $t=0$),
 and the condition
$$
M_{K_1}(u)=M_{K_2}(u)\qquad \forall u\in {\mathbb S}^{d-1}
$$
implies $K_1=K_2$.

It is also well known, \cite{BF},  that on the plane there are convex bodies $K$ that are {\it not} Euclidean discs, but nevertheless  satisfy $M_K(u)=1$  for all $u\in {\mathbb S}^1$.  These are the {\it bodies of  constant width} $1$.

In  1969 V. Klee asked whether the condition $M_{K_1}\equiv M_{K_2}$ implies that $K_1$ and $K_2$ are essentially the same (i.e.,  differ by translation and/or reflection about the origin)
in general, or, at least, whether the condition $M_K\equiv c$ implies that  $K$ is a Euclidean ball, \cite{Kl}.

Recently,  R. Gardner and V. Yaskin, together with  the second and the third named authors gave the negative  answer to the first question of V. Klee by
 constructing two  bodies of revolution $K_1$, $K_2$ such that $K_1$ is origin-symmetric, $K_2$ is not origin-symmetric, but $M_{K_1}\equiv M_{K_2}$  (see \cite{GRYZ}). In \cite{NRZ} this result was  strengthened  for all  {\it even} dimensions. It was shown that there exist two essentially different convex bodies of revolution $K_1$, $K_2\subset {\mathbb R}^d$ such that $A_{K_1}\equiv A_{K_2}$, $M_{K_1}\equiv M_{K_2}$, and $P_{K_1}\equiv P_{K_2}$, where, for $u\in {\mathbb S}^{d-1}$,
$$
A_K(u)=\mathrm{vol}_{d-1}(K\cap u^{\perp}),\qquad P_K(u)=\mathrm{vol}_{d-1}(K|u^{\perp}),
$$
and
$K|u^{\perp}$ is the projection of $K$ to $u^{\perp}$.

In this paper we answer the second question of V. Klee. Our main result is the following

\bt\label{klee}
If  $d\ge 3$, there exists a convex body of revolution $K\subset {\mathbb R}^d$ satisfying $M_K\equiv const$ that is not a Euclidean ball.
\et
\begin{remark}\label{hryusha}
An alternative proof in  the case $d=4$ has been  given in \cite{NRZ}. Unfortunately, the elementary techniques used there fail for other dimensions.
\end{remark}
\medskip

Our bodies  will be {\it small perturbations of the Euclidean ball}.
The proofs of Theorem~\ref{klee} for even and odd $d$ are different.  For even $d$  we can get away with elementary algebra and
  control finitely many  moments of $f$  to obtain a  {\it local} perturbation of the ball.  The case $d\ge 3$, $d$ is odd,  is much more involved.
 To control  the perturbation,
 we use the properties of the  {\it Spherical Radon Transform},
together with
the  {\it Borsuk-Ulam Theorem}  asserting that any continuous map from ${\mathbb S}^m$ to ${\mathbb R}^m$,  taking  antipodal points to antipodal points, contains zero in its image. The reader can find all  necessary  information in \cite{He} and  \cite{Mat}.

The paper is organized as follows. In  Section 2 we  reduce the problem to finding a non-trivial solution to two integral equations.
 In Section 3   we prove Theorem \ref{klee} for even  $d$.
Section 4 is devoted to the most difficult part of the proof when $d$ odd.
 The Appendix  contains technical parts of the proofs and some auxiliary statements.

\section{Reduction to a system of integral equations}

From now on,  we assume that $d\ge 3$. We will be dealing with convex    bodies of revolution
$$
K_f=\{x\in \R^d:\,x_1\in[-\lambda,\mu],\,\,x_2^2+x_3^2+...+x_d^2\le f^2(x_1)\},
$$
obtained by the rotation of   a smooth (except for the endpoints) strictly concave function $f:[-\lambda,\mu]\to [0,\infty)$ about the $x_1$-axis, where $\lambda$ and $\mu$  are some positive real numbers.

Note that $K$ is rotation invariant, thus every its hyperplane section is
equivalent to a section by a hyperplane with normal vector in the second quadrant of the $(x_1,x_2)$-plane.

\pagebreak

\begin{figure}[ht]
\includegraphics[width=420pt]{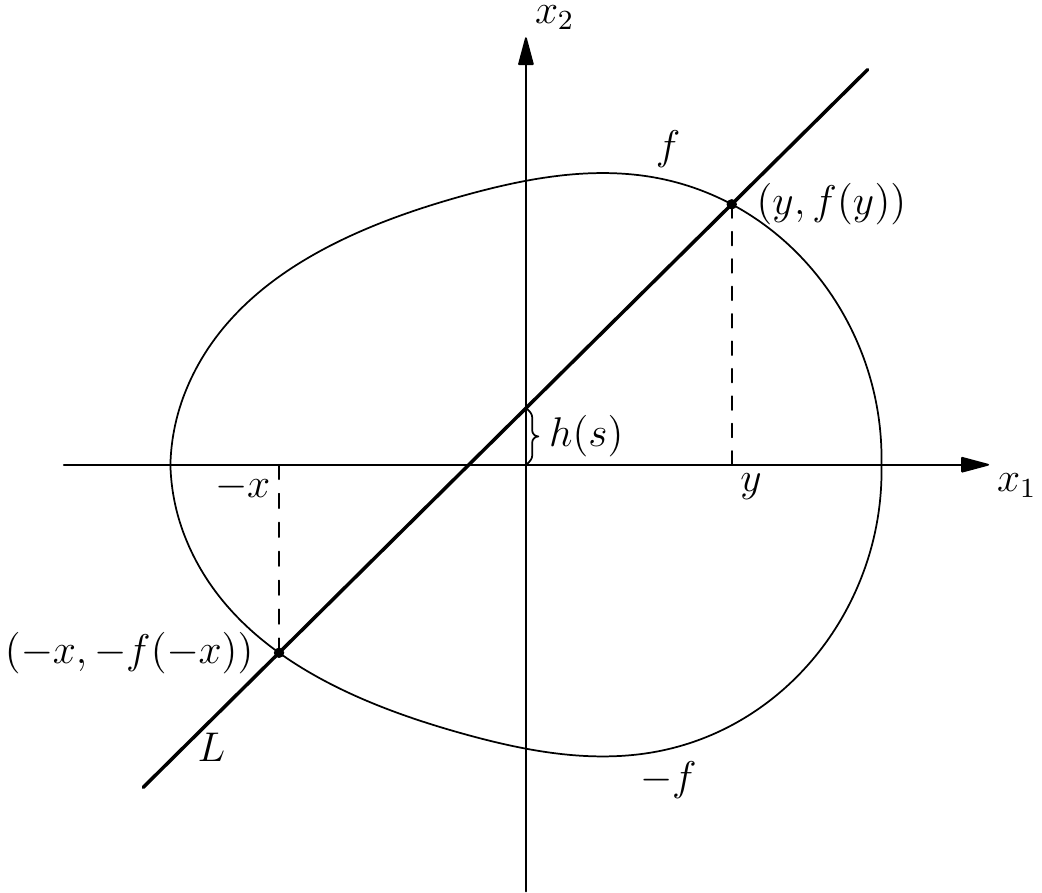}
\caption{Sections of $K$ and $H(L)$  by the  $(x_1,x_2)$-plane.}
\label{pic0}
\end{figure}

\bl\label{fn}
 Let $L(\xi)=L(s,h,\xi)=s\xi+h$ be a linear function with slope $s$, and let
  $H(L)=\{x\in{\mathbb R}^d:\,x_2=L(x_1)\}$ be the corresponding hyperplane.
  Then the   section $K\cap H(L)$ is of maximal volume if and only if
  \begin{equation}\label{max}
\int\limits_{-x}^{y}(f^2-L^2)^{(d-4)/2}L=0,
\end{equation}
where $-x$ and $y$ are the first coordinates of the points at which $L$ intersects the graphs of $-f$ and $f$
respectively (see Figure \ref{pic0}).
\el

\bp
Fix $s>0$. Observe  that
 the  slice $K\cap H(L)\cap H_{\xi}$ of $K\cap H(L)$ by the hyperplane $H_{\xi}=\{x\in {\mathbb R}^d:\,x_1=\xi\}$, $-x(s)<\xi<y(s)$,    is the
$(d-2)$-dimensional Euclidean ball $\{(\xi,L(\xi),x_3,x_4,...,x_d):\,x_3^2+...+x^2_d\le r^2\}$ of radius $r=\sqrt{f^2(\xi)-L^2(\xi)}$.
Hence,
\begin{equation}\label{const1}
\mathrm{vol}_{d-1}(K \cap H(L))=v_{d-2}\sqrt{1+s^2} \int\limits_{-x(s)}^{y(s)}(f^2(\xi)-L^2(\xi))^{(d-2)/2}d\xi,
\end{equation}
where $v_{d-2}$ is the volume of the unit ball in $\R^{d-2}$.

The section $K\cap H(L)$ is  of maximal volume if and only if
$$
\frac{d}{d h}\mathrm{vol}_{d-1}(K \cap H(L))=0,
$$
where in  the {\it only if} part we use the Theorem of  Brunn, \cite{Ga}. Computing the derivative, we conclude that for a  given $s\in {\mathbb R}$, the section $K\cap H(L)$ is  of maximal volume if and only if (\ref{max}) holds.
\ep
\bl\label{fn1}
 Let $L(s,\xi)=s\xi+h(s)$ be a family of linear functions parameterized by the slope $s$. For each $L$ in our family, define the hyperplane $H(L)$ by
  $H(L)=\{x\in{\mathbb R}^d:\,x_2=L(x_1)\}$ (see Figure \ref{pic0}). The corresponding family of
 sections is of constant $(d-1)$-dimensional volume  if and only if
 \begin{equation}\label{const}
\int\limits_{-x}^{y}(f^2-L^2)^{(d-2)/2}=\frac{const}{\sqrt{1+s^2}} \qquad \textrm{for all}\qquad s>0.
\end{equation}
In the case of the unit Euclidean ball,  the constant is equal to  $\frac{v_{d-1}}{v_{d-2}}$.
\el
\bp
The right hand side in (\ref{const1}) is constant if and only if  (\ref{const}) holds.
\ep

In what follows we will choose  the {\it perturbation function} $h$ that is  infinitely smooth,  not identically zero, supported on $[1-2\delta, 1-\delta]$ for some small $\delta>0$, and is small together with sufficiently many derivatives. The convex body corresponding to any such function will be automatically  asymmetric since not all its maximal sections will pass through a single point.

\section{The case  of even $d$ }

Note that in this case $p=\frac{d-2}{2}\in{\mathbb N}$.
Then (\ref{const}) and (\ref{max}) take the form
\begin{equation}\label{devencon}
\int\limits_{-x}^{y}(f^2-L^2)^p=\int\limits_{-x_o}^{y_o}(f_o^2-L_o^2)^p=\frac{const}{\sqrt{1+s^2}},
\end{equation}
\begin{equation}\label{devenmax}
   \int\limits_{-x}^{y}(f^2-L^2)^{p-1}L =0,
\end{equation}
where $f_o(\xi)=\sqrt{1-\xi^2}$, $L_o(s,\xi)=s\xi$, and  $y_o(s)=x_o(s)=1/\sqrt{1+s^2}$.

\pagebreak

\begin{figure}[ht]
\includegraphics[width=420pt]{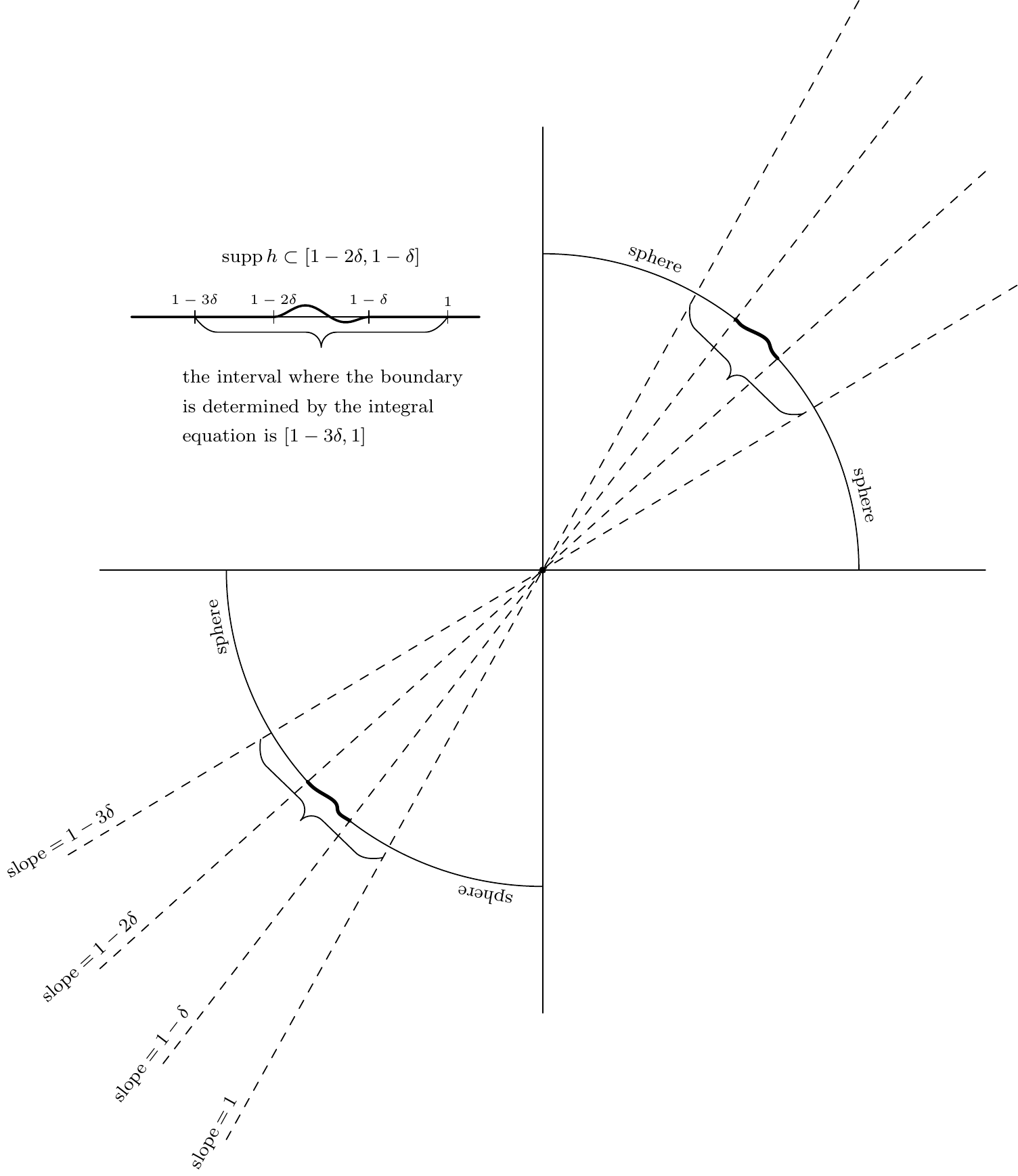}
\caption{Graph of  $f$ (the case of an even dimension).}
\label{pic1}
\end{figure}

Our body of revolution $K_f$ will be constructed as a {\it local} perturbation of the Euclidean ball (see Figure \ref{pic1}).

\pagebreak

\begin{figure}[ht]
\includegraphics[width=350pt]{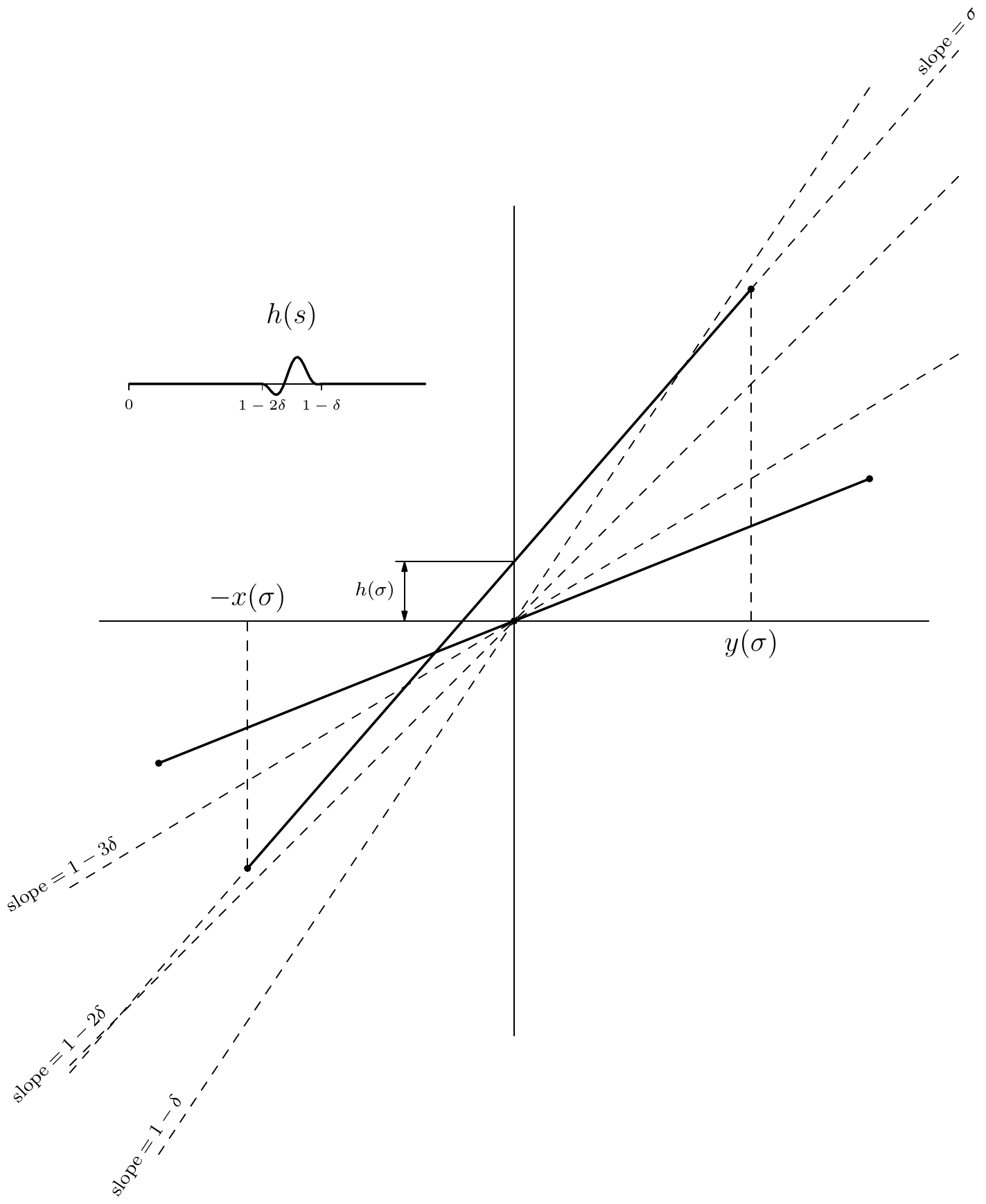}
\caption{The functions $x(\sigma)$ and $y(\sigma)$.}
\label{pic4}
\end{figure}

The equations
\begin{equation}\label{ura}
f(y(\sigma))=L(\sigma,y(\sigma)),\qquad f(-x(\sigma))=L(\sigma,-x(\sigma))
\end{equation}
show that to define $f$, it is enough to define two decreasing functions $x(\sigma)$ and $y(\sigma)$ on $[0,+\infty)$ (see Figure \ref{pic4}).

\setstretch{1.22}

For the unperturbed case of the unit ball,  $h\equiv 0$ and these functions are just $y_o(s)=x_o(s)=1/\sqrt{1+s^2}$. Our new functions $x(\sigma)$ and $y(\sigma)$ will coincide with $x_o$ and $y_o$ for all $\sigma\notin[1-2\delta,1-\delta]$. Since the curvature of the semicircle is strictly positive, the resulting function $f$ will be strictly concave if $x$ and $y$ are close to $x_o$ and $y_o$ in $C^2$.

We shall make our construction in several steps. First, we {\it define} $x=x_o$, $y=y_o$ on $[1,\infty)$. Second, we will transfrom   equations (\ref{devencon}), (\ref{devenmax})  to obtain equations (\ref{hah1}), (\ref{hah2}) written purely in terms of $x$ and $y$ (see below). Then we
will use these new equations
 to extend the functions $x$ and $y$ to $[1-3\delta,1]$. Due to the existence and uniqueness lemma  and the  remark after it (Lemma \ref{dif} and Remark \ref{r} in Appendix), we will be able to do it if $\delta$ and $h$ are  sufficiently small, and, moreover, the extensions will coincide with $x_o$ and $y_o$ on $[1-\delta,1]$ and will be close to $x_o$ and $y_o$ with two derivatives on $[1-3\delta,1-\delta]$. Then, we will show that a little miracle happens and our extensions automatically coincide with $x_o$ and $y_o$ on $[1-3\delta,1-2\delta]$ as well. This will allow us to put $x=x_o$, $y=y_o$ on the remaining interval $[0, 1-3\delta]$ and get a nice smooth function (see Figure \ref{pic1}). At last, we will show that equations (\ref{devencon}), (\ref{devenmax}) will hold up to $s=0$, thus finishing the proof.

\medskip

{\bf Step 1}.
We put   $x=x_o$, $y=y_o$ on $[1,\infty)$.

\medskip

{\bf Step 2}. To construct $x$, $y$ on $[1-3\delta,1]$,   we will make  some technical preparations.
First,  we will  differentiate equations (\ref{devencon}), (\ref{devenmax})  a few times
to obtain a   system of {\it four} integral  equations with  {\it four} unknown functions $x$, $y$, $x'$, $y'$.
  Next, we will apply Lemma \ref{dif} and Remark \ref{r} to show
 that there exists a  solution $x$,  $y$, $x'$, $y'$ of the constructed system of integral equations on $[1-3\delta,1]$,  which coincides with $x_o$, $y_o$, $\frac{dx_o}{ds}$, $\frac{dy_o}{ds}$ on $[1-\delta,1]$.
 Finally, we will prove that the $x$ and $y$ components of that solution  give a solution of (\ref{devencon}), (\ref{devenmax})  with $f$ defined by (\ref{ura}).

Differentiating  equation  (\ref{devencon}) $p+1$ times and   equation   (\ref{devenmax}) $p$ times,
 we obtain
$$
(-2)^pp!\Big[\Big(\Big(L\frac{\partial L}{\partial s}\Big)\Big|_{(s,-x(s))}\Big)^p\frac{dx}{ds}(s)\,+
\Big(\Big(L\frac{\partial L}{\partial s}\Big)\Big|_{(s,y(s))}\Big)^p\frac{dy}{ds}(s)\Big]\,+
$$
\begin{equation}\label{feq}
\int\limits_{-x(s)}^{y(s)}\Big(\frac{\partial}{\partial s}\Big)^{p+1}\Big((f^2(\xi)-L^2(s, \xi))^p\Big)d\xi\,=\Big(\frac{d}{ds}\Big)^{p+1}\Big(\frac{const}{\sqrt{1+s^2}}\Big),
\end{equation}
and
$$
(-2)^{p-1}(p-1)!\Big[\Big(\Big(L\frac{\partial L}{\partial s}\Big)^{p-1}L\Big)\Big|_{(s,-x(s))}\frac{dx}{ds}(s)\,+
\Big(\Big(L\frac{\partial L}{\partial s}\Big)^{p-1}L\Big)\Big|_{(s,y(s))}\frac{dy}{ds}(s)\Big]\,+
$$
\begin{equation}\label{seq}
\int\limits_{-x(s)}^{y(s)}\Big(\frac{\partial}{\partial s}\Big)^{p}\Big((f^2(\xi)-L^2(s, \xi))^{p-1}L(s,\xi)\Big)d\xi=0.
\end{equation}
When $s\le 1$, the integral term $I$  in  (\ref{feq}) can be split as
$$
I=\int\limits_{-x(s)}^{y(s)}\Big(\frac{\partial}{\partial s}\Big)^{p+1}\Big((f^2(\xi)-L^2(s, \xi))^p\Big)d\xi=
$$
$$
\Big(\int\limits_{-x(s)}^{-x_o(1)}\,+\,\int\limits_{y_o(1)}^{y(s)}\Big)\,\Big(\frac{\partial}{\partial s}\Big)^{p+1}\Big((f^2(\xi)-L^2(s, \xi))^p               \Big)d\xi\,+\,\Xi_1(s),
$$
where
$$
\Xi_1(s)=\int\limits^{y_o(1)}_{-x_o(1)}\Big(\frac{\partial}{\partial s}\Big)^{p+1}\Big((f_o^2(\xi)-L^2(s, \xi))^p\Big)d\xi.
$$
Making the change of variables $\xi=-x(\sigma)$ in the integral $\int_{-x(s)}^{-x_o(1)}$, and $\xi=y(\sigma)$  in the integral $\int_{y_o(1)}^{y(s)}$, we obtain
$$
I=-\int\limits_s^{1}\Big(\frac{\partial}{\partial s}\Big)^{p+1}\Big(L^2(\sigma,-x(\sigma))-L^2(s,-x(\sigma))\Big)^p\,\frac{dx}{ds}(\sigma)d\sigma\,-
$$
$$
\int\limits_s^{1}\Big(\frac{\partial}{\partial s}\Big)^{p+1}\Big(L^2(\sigma,y(\sigma))-L^2(s,y(\sigma))\Big)^p\,\frac{dy}{ds}(\sigma)d\sigma\,+\,\Xi_1(s).
$$
Similarly,  we have
$$
\int\limits_{-x(s)}^{y(s)}\Big(\frac{\partial}{\partial s}\Big)^{p}\Big((f^2(\xi)-L^2(s, \xi))^pL(s,\xi)\Big)d\xi=
$$
$$
-\int\limits_s^{1}\Big(\frac{\partial}{\partial s}\Big)^{p}\Big(\Big(L^2(\sigma,-x(\sigma))-L^2(s,-x(\sigma))\Big)^{p-1}L(s,-x(\sigma))\Big)\,\frac{dx}{ds}(\sigma)d\sigma\,-
$$
$$
\int\limits_s^{1}\Big(\frac{\partial}{\partial s}\Big)^{p}\Big(\Big(L^2(\sigma,y(\sigma))-L^2(s,y(\sigma))\Big)^{p-1}L(s,y(\sigma))\Big)\,\frac{dy}{ds}(\sigma)d\sigma\,+\,\Xi_2(s),
$$
where
$$
\Xi_2(s)=\int\limits^{y_o(1)}_{-x_o(1)}\Big( \frac{\partial}{\partial s}\Big)^{p}\Big((f_o^2(\xi)-L^2(s, \xi))^{p-1}L(s,\xi)                 \Big)d\xi.
$$
To reduce the resulting system of integro-differential equations to a pure system of integral equations we add two independent unknown functions $x'$, $y'$ and
two new relations
$$
x(s)=-\int\limits_s^{1}x'(\sigma)d\sigma +x_o(1),\qquad y(s)=-\int\limits_s^{1}y'(\sigma)d\sigma +y_o(1).
$$
We rewrite   our  equations (\ref{devencon}), (\ref{devenmax}) as follows:
\begin{equation}\label{hah1}
(-2)^pp!\Big[\Big(\Big(L\frac{\partial L}{\partial s}\Big)\Big|_{(s,-x(s))}\Big)^px'(s)\,+
\Big(\Big(L\frac{\partial L}{\partial s}\Big)\Big|_{(s,y(s))}\Big)^py'(s)\Big]\,-
\end{equation}
$$
\int\limits_s^{1}\Big(\frac{\partial}{\partial s}\Big)^{p+1}\Big(L^2(\sigma,-x(\sigma))-L^2(s,-x(\sigma))\Big)^p\,x'(\sigma)d\sigma\,
-
$$
$$
\int\limits_s^{1}\Big(\frac{\partial}{\partial s}\Big)^{p+1}\Big(L^2(\sigma,y(\sigma))-L^2(s,y(\sigma))\Big)^p\,y'(\sigma)d\sigma\,+\,\Xi_1(s)=\Big(\frac{d}{ds}\Big)^{p+1}\Big(\frac{const}{\sqrt{1+s^2}}\Big),
$$
and
\begin{equation}\label{hah2}
(-2)^{p-1}(p-1)!\Big[\Big(\Big(L\frac{\partial L}{\partial s}\Big)^{p-1}L\Big)\Big|_{(s,-x(s))}x'(s)\,+
\Big(\Big(L\frac{\partial L}{\partial s}\Big)^{p-1}L\Big)\Big|_{(s,y(s))}y'(s)\Big]
\,-
\end{equation}
$$
\int\limits_s^{1}\Big(\frac{\partial}{\partial s}\Big)^{p}\Big(\Big(L^2(\sigma,-x(\sigma))-L^2(s,-x(\sigma))\Big)^{p-1}L(s,-x(\sigma))\Big)\,x'(\sigma)d\sigma\,-
$$
$$
\int\limits_s^{1}\Big(\frac{\partial }{\partial s}\Big)^{p}\Big(\Big(L^2(\sigma,y(\sigma))-L^2(s,y(\sigma))\Big)^{p-1}L(s,y(\sigma))\Big)\,y'(\sigma)d\sigma\,+\Xi_2(s)=0.
$$
Now we  rewrite our system  in the  form
\begin{equation}\label{hah41}
\mathbf{G}(s,Z(s))=\int\limits_s^{1}\mathbf{\Theta}(s,\sigma, Z(\sigma))d\sigma+\mathbf{\Xi}(s).
\end{equation}
Here
$$
Z=\left(
\begin{array}{cccc}
x\\
y\\
x'\\
y'
\end{array}
\right),
$$
$$
\mathbf{G}(s,Z)=\left(
\begin{array}{cccc}
x\\
y\\[7pt]
(-2)^pp!\Big[\Big(L\frac{\partial L}{\partial s}\Big|_{(s,-x)}\Big)^px'\,+
\Big(L\frac{\partial L}{\partial s}\Big|_{(s,y)}\Big)^py'\Big]\\[15pt]
(-2)^{p-1}(p-1)!\Big[\Big(\Big(L\frac{\partial L}{\partial s}\Big)^{p-1}L\Big)\Big|_{(s,-x)}x'\,+
\Big(\Big(L\frac{\partial L}{\partial s}\Big)^{p-1}L\Big)\Big|_{(s,y)}y'\Big]\\
\end{array}
\right),
$$
$$
\mathbf{\Theta}(s,\sigma, Z)=
-\left(
\begin{array}{cccc}
x'\\
y'\\
\Theta_{1}\\
\Theta_{2}
\end{array}
\right),
$$
where
$$
\Theta_{1}=-\Big(\frac{\partial}{\partial s}\Big)^{p+1}\Big(L^2(\sigma,-x)-L^2(s,-x)\Big)^p\,x'\,
-
\Big(\frac{\partial}{\partial s}\Big)^{p+1}\Big(L^2(\sigma,y)-L^2(s,y)\Big)^p\,y'\,,
$$
$$
\Theta_{2}=
-\Big(\frac{\partial}{\partial s}\Big)^{p}\Big(\Big(L^2(\sigma,-x)-L^2(s,-x)\Big)^{p-1}L(s,-x)\Big)\,x'\,-
$$
$$
\Big(\frac{\partial}{\partial s}\Big)^{p}\Big(\Big(L^2(\sigma,y)-L^2(s,y)\Big)^{p-1}L(s,y)\Big)\,y',
$$
and
$$
\mathbf{\Xi}(s)=\left(
\begin{array}{cccc}
x_o(1)\\
y_o(1)\\
-\Xi_1(s)+\Big(\frac{d}{ds}\Big)^{p+1}\Big(\frac{const}{\sqrt{1+s^2}}\Big)\\
-\Xi_2(s)
\end{array}
\right).
$$
 Note that ${\mathbf G}$, ${\mathbf \Theta}$, ${\mathbf \Xi}$ are well-defined and infinitely smooth for all $s,\sigma \in(0,1]$ and $Z\in \R^4$.
Observe also that
$$
D_Z{\mathbf G}\Big|_{(s,Z)}=
\left(
\begin{array}{cc}
\mathbf{I}&0\\
\mathbf{*}&\mathbf{A}\\
\end{array}
\right),
$$
where
$$
{\mathbf I}=\left(
\begin{array}{ccc}
1 & 0\\
0&
1
\end{array}
\right),\qquad{\mathbf A}={\mathbf A(s,x,y)}=
$$
$$
\left(
\begin{array}{ccc}
(-2)^pp!\Big(\Big(L\frac{\partial L}{\partial s}\Big)\Big|_{(s,-x)}\Big)^p\qquad\qquad & (-2)^pp!\Big(\Big(L\frac{\partial L}{\partial s}\Big)\Big|_{(s,y)}\Big)^p\\
(-2)^{p-1}(p-1)!\Big(\Big(L\frac{\partial L}{\partial s}\Big)^{p-1}L \Big)\Big|_{(s,-x)}\qquad\qquad&
(-2)^{p-1}(p-1)!\Big(\Big(L\frac{\partial L}{\partial s}\Big)^{p-1}L\Big)\Big|_{(s,y)}
\end{array}
\right).
$$

\smallskip

The function
$$
Z_o(s)=\left(
\begin{array}{cccc}
x_o(s)\\
y_o(s)\\[3pt]
\frac{dx_o}{ds}(s)\\[3pt]
\frac{dy_o}{ds}(s)
\end{array}
\right)
$$
solves the system (\ref{hah41}) with ${\mathbf G}$, ${\mathbf \Theta}$, ${\mathbf \Xi}$ corresponding to $h\equiv 0$ (we will denote them by
${\mathbf G_o}$, ${\mathbf \Theta_o}$, ${\mathbf \Xi_o}$) on $[\frac{1}{2},1]$, say.

We claim that
\begin{equation}\label{detH}
\det\,\Big(D_Z{\mathbf G_o}\Big|_{(s,Z_o(s))}  \Big)=\det({\mathbf A_o(s,x_o(s),y_o(s))})\neq 0\qquad \textrm{for all}\,s\in(0,1].
\end{equation}
Indeed, the  matrix ${\mathbf A_o}(s,x_o(s),y_o(s))$
has the sign pattern
$$
\left(
\begin{array}{ccc}
+ & +\\
+ &  -
\end{array}
\right),
  \qquad \textrm{when}\,\, p\,\,\textrm{ is even},\quad \textrm{and }
\qquad
\left(\begin{array}{ccc}
- & -\\
-&
+
\end{array}
\right),\quad \textrm{when}\, p\,\,\textrm{ is odd}.
$$
Thus, (\ref{detH}) follows.
In particular,
$$
\det\,\Big(D_Z{\mathbf G_o}\Big|_{(1,Z_o(1))}  \Big)\neq 0.
$$

Lemma \ref{dif} implies then that we can choose some small $\delta>0$ and, for any fixed $k\in {\mathbb N}$,  construct a $C^k$-close to $Z_o(s)$ solution $Z(s)$ of (\ref{hah41}) on $[1-3\delta,1]$ whenever
${\mathbf G}$, ${\mathbf \Theta}$, ${\mathbf \Xi}$ are sufficiently close to ${\mathbf G_o}$, ${\mathbf \Theta_o}$, ${\mathbf \Xi_o}$ in $C^k$ on  certain compact sets. Since
$\mathbf{G}$, $\mathbf{\Theta}$, $\mathbf{\Xi}$ and their derivatives are some explicit (integrals of) polynomials in $Z$, $s$, $\sigma$, $h(s)$, and the derivatives of $h(s)$,
this closeness assumption will hold if $h$ is sufficiently close to zero with sufficiently many derivatives.
Moreover, since $h$ vanishes on $[1-\delta,1]$, the assumptions of Remark \ref{r} are satisfied and we have $Z(s)=Z_o(s)$ on $[1-\delta,1]$.

To prove that the $x$ and $y$ components of the solution we found  give a solution of (\ref{devencon}), (\ref{devenmax})  with $f$ defined by (\ref{ura}), we
consider the functions
$$
F(s):=
\int\limits_{-x(s)}^{y(s)}\Big(f(s,\xi)^2-L^2(s,\xi)\Big)^pd\xi-\frac{const}{\sqrt{1+s^2}},
$$
$$
 H(s):=  \int\limits_{-x(s)}^{y(s)}\Big(f(s,\xi)^2-L(s,\xi)^2\Big)^{p-1}L(s,\xi)d\xi.
$$
Since equations (\ref{hah1}) and  (\ref{hah2}) of our system (\ref{hah41}) were obtained by the  differentiation of equations (\ref{devencon}),  (\ref{devenmax}),  we have
$$
\Big(\frac{d}{ds}\Big)^{p+1}F(s)=0,\qquad \Big(\frac{d}{ds}\Big)^pH(s)=0
$$
on $[1-3\delta, 1]$. Hence, $F$ and $H$ are polynomials
on $[1-3\delta, 1]$.
Since $h(s)=0$, $x(s)=x_o(s)$, $y(s)=y_o(s)$  on  $[1-\delta, 1]$,  $F$ and $H$  vanish on $[1-\delta, 1]$ and, therefore, identically.
Thus, we conclude
 that  the $x$ and $y$ components of    the solutions
of (\ref{hah1}), (\ref{hah2}) solve  (\ref{devencon}),  (\ref{devenmax})  on $(1-3\delta,1]$.
 Step 2 is completed.

\medskip

{\bf Step 3}.
 We claim that $x=x_o$, $y=y_o$ on $[1-3\delta, 1-2\delta]$, i.e., the {\it perturbed solution} returns to the semicircle.
Since  $h$ is supported on $[1-2\delta,1-\delta]$,  we have $L=L_o=s\xi$ and
$\frac{\partial}{\partial s}L(s,\xi)=\xi$ for $s\in [1-3\delta,1-2\delta]$. It follows that every time we differentiate equation   (\ref{devencon})  (with respect to $s$)
 we can divide the result by $s$ to obtain
\begin{equation}\label{evenmom}
\int\limits_{-x(s)}^{y(s)}(f^2(\xi)-L_o^2(s,\xi))^{p-k}\xi^{2k}d\xi=
\int\limits_{-x_o(s)}^{y_o(s)}(f_o^2(\xi)-L_o^2(s,\xi))^{p-k}\xi^{2k}d\xi,\qquad k\le p.
\end{equation}
If we take $k=p$ in (\ref{evenmom}), we get
\begin{equation}\label{evenk1}
\int\limits_{-x(s)}^{y(s)}\xi^{2p}d\xi=
\int\limits_{-x_o(s)}^{y_o(s)}\xi^{2p}d\xi.
\end{equation}
Similarly, for $k\le p-1$,   equation (\ref{devenmax}) implies that
\begin{equation}\label{oddmom}
\int\limits_{-x(s)}^{y(s)}(f^2(\xi)-L_o^2(s,\xi))^{p-1-k}\xi^{2k+1}d\xi=0.
\end{equation}
 Putting $k=p-1$ in (\ref{oddmom}), we get
\begin{equation}\label{oddk}
\int\limits_{-x(s)}^{y(s)}\xi^{2p-1}d\xi=0=
\int\limits_{-x_o(s)}^{y_o(s)}\xi^{2p-1}d\xi.
\end{equation}
Equation (\ref{oddk}) yields $x(s)=y(s)$, and the symmetry (with respect to $0$) of intervals $(-x_o(s),y_o(s))$, $(-x(s),y(s))$,
together with (\ref{evenk1}), yield $(-x_o(s),y_o(s))=(-x(s),y(s))$ for all $s\in [1-3\delta, 1-2\delta]$.
Step 3 is completed.

\medskip

{\bf Step 4}.  We put $x=x_o$, $y=y_o$ on $ [0,1-3\delta]$, which will result in a function $f$ defined  on $[-1,1]$ and coinciding with $\sqrt{1-\xi^2}$ outside
small intervals around $\pm\frac{1}{\sqrt{2}}$.
It remains to check that (\ref{devencon}),  (\ref{devenmax}) are valid for $s\in[0,1-3\delta]$.  We  will prove the validity of  (\ref{devencon}). The proof
for  equation (\ref{devenmax})  is similar.

Since $h\equiv 0$ away from $(1-2\delta,1-\delta)$, we have $L(s,\xi)=s\xi$ for $s\in[0,1-3\delta]$,  so we need to check that
$$
\int\limits_{-x(s)}^{y(s)}(f^2(\xi)-(s\xi)^2)^pd\xi=\int\limits_{-x(s)}^{y(s)}(f_o^2(\xi)-(s\xi)^2)^pd\xi,\qquad \forall s\in[0,1-3\delta].
$$
Recall that $x=x_o$ and $y=y_o$ everywhere on this interval, so we can write $x$ and $y$ instead of $x_o$ and $y_o$ on the right hand side.

Opening the parentheses, we see that it suffices to check that
\begin{equation}\label{devenm}
\int\limits_{-x(s)}^{y(s)}f^{2j}(\xi)\xi^{2(p-j)}d\xi=\int\limits_{-x(s)}^{y(s)}f_o^{2j}(\xi)\xi^{2(p-j)}d\xi,\quad \forall j=1,\dots,p,\quad s\in[0,1-3\delta].
\end{equation}
Since $f\equiv f_o$  outside $[-x(1-3\delta),y(1-3\delta)]$, it is enough to check (\ref{devenm}) for
$s=1-3\delta$.

 To this end, we first  take $s=1-3\delta$, $k=p-1$ in (\ref{evenmom}) and conclude that
\begin{equation}\label{mom2}
\int\limits_{-x(1-3\delta)}^{y(1-3\delta)}f^{2}(\xi)\xi^{2p-2}d\xi=\int\limits_{-x(1-3\delta)}^{y(1-3\delta)}f_o^{2}(\xi)\xi^{2p-2}d\xi,
\end{equation}
which is (\ref{devenm}) for $j=1$.
Now we go ``one step up",  by taking $s=1-3\delta$, $k=p-2$ in (\ref{evenmom}), to get
$$
\int\limits_{-x(1-3\delta)}^{y(1-3\delta)}(f^{2}(\xi)-(s\xi)^2)^2\xi^{2p-4}d\xi=\int\limits_{-x(1-3\delta)}^{y(1-3\delta)}(f_o^{2}(\xi)-(s\xi)^2)^2\xi^{2p-4}d\xi.
$$
The last equality  together with (\ref{mom2}) yield
$$
\int\limits_{-x(1-3\delta)}^{y(1-3\delta)}f^{4}(\xi)\xi^{2p-4}d\xi=\int\limits_{-x(1-3\delta)}^{y(1-3\delta)}f_o^{4}(\xi)\xi^{2p-4}d\xi,
$$
which is (\ref{devenm}) for $j=2$. Proceeding in a similar way we get (\ref{devenm}) for $j=1,\dots,p$.
This finishes the proof of Theorem \ref{klee} in even dimensions.

\section{The  odd dimensional case}

Note that in this case $p=q+\frac{1}{2}$, $q\in{\mathbb Z_+}$.
Then (\ref{const}) and (\ref{max}) take the form
\begin{equation}\label{dodd+}
\int\limits_{-x}^{y}(f^2-L^2)^{q+\frac{1}{2}}=\int\limits_{-x_o}^{y_o}(f_o^2-L_o^2)^{q+\frac{1}{2}}=\frac{const}{\sqrt{1+s^2}},
\end{equation}
\begin{equation}\label{dodd-}
  \int\limits_{-x}^{y}(f^2-L^2)^{q-\frac{1}{2}}L =0,
\end{equation}
where $f_o(\xi)=\sqrt{1-\xi^2}$, $L_o(s,\xi)=s\xi$, and  $y_o(s)=x_o(s)=1/\sqrt{1+s^2}$.

Let $L=L(s, \xi)=s \xi+h(s)$ be  a family of  linear functions, parameterized by the slope $s$.
Here the {\it perturbation function} $h$ is  infinitely smooth, supported on $[1-2\delta, 1-\delta]$ for some small $\delta>0$ to be chosen later, and is small together with sufficiently many derivatives.

\pagebreak

\begin{figure}[ht]
\includegraphics[width=350pt]{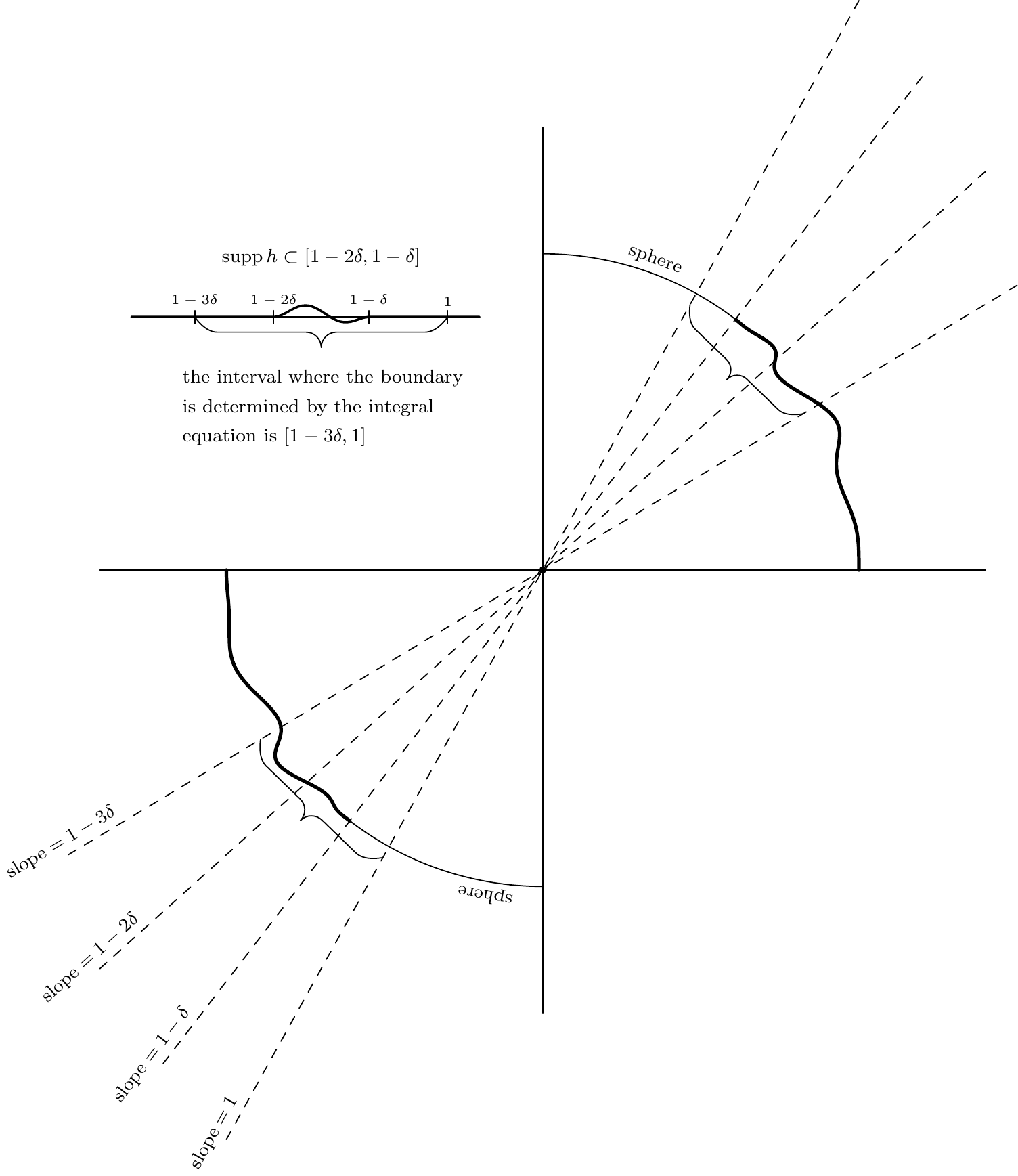}
\caption{Graph of  $f$ (the case of an odd dimension).}
\label{pic5}
\end{figure}

Our body of revolution $K_f$ will be constructed as a  perturbation of the Euclidean ball (see Figure \ref{pic5},  and compare it with Figure  \ref{pic1}). Note that in the case of an odd dimension,  the perturbation is not local.

\setstretch{1.02}

We shall make our construction in several steps  corresponding to the slope ranges $s\in [1,\infty)$, $s\in[1-3\delta,1]$, and $s\in (0,1-3\delta]$.
We will use different   ways to describe the boundary of $K_f$ within those ranges. We will {\it define}
$f(\xi)=f_o(\xi)$ for  $\xi\in\left[-\frac{1}{\sqrt{2}},\frac{1}{\sqrt{2}}\right]$.
We will differentiate    (\ref{dodd+}), (\ref{dodd-}) and rewrite  the resulting equations in terms of $x$ and $y$ to extend
$x$ and $y$ to $[1-3\delta, 1]$ like we did in the even case. As before, $f$ is related to $x$ and $y$ by (\ref{ura}).
Finally, we will change the point of view and define the remaining part of $f$ in terms of the radial functions
$R(\alpha)$ and $r(\alpha)$, related to $f$  by
\begin{equation}\label{f}
f(R(\alpha)\cos\alpha)=R(\alpha)\sin\alpha,\qquad f(-r(\alpha)\cos\alpha)=r(\alpha)\sin\alpha,\qquad \alpha\in[0,\tfrac{\pi}{2}].
\end{equation}
Note that the {\it radial function }
$$
\rho_K(u)=\sup\{t>0:\,\,tu\in K\}
$$
of the resulting body $K$ satisfies
\begin{equation}\label{tryuk}
\rho_K(u)= \begin{cases}
R(\alpha)\qquad &\text{if \,}u_1>0,\\
r(\alpha) & \text{if \,}u_1<0,
\end{cases}
\end{equation}
where $u=(u_1,\dots)\in {\mathbb S}^{d-1}$ and $\alpha\in[0,\frac{\pi}{2}]$, $\cos\alpha=|u_1|$.

The solutions $R(\alpha)$ and $r(\alpha)$ of the equations that we will use during the last step may develop a singularity at $\alpha=0$.
To avoid this, we will impose several additional cancellation conditions on
 the perturbation function $h$. We will  use the Borsuk-Ulam Theorem    to show the existence of a non-zero function $h$ satisfying these cancellation restrictions.

\medskip

{\bf Step 1}.
We put $x=x_o$, $y=y_o$ on $[1,\infty)$, which is equivalent to putting $f(\xi)=\sqrt{1-\xi^2}$ for $\xi\in[-\frac{1}{\sqrt{2}},\frac{1}{\sqrt{2}}]$.

\medskip

{\bf Step 2}.
Differentiating  equation  (\ref{dodd+}) $q+1$ times, we obtain
\begin{equation}\label{ass+}
\Big(\frac{\partial}{\partial s}  \Big)^{q+1}\int\limits_{-x(s)}^{y(s)}(f^2(\xi)-L^2(s,\xi))^{q+\frac{1}{2}}d\xi=
\end{equation}
$$
\Big(\int\limits_{-x(s)}^{-x_o(1)}\,+\,\int\limits_{y_o(1)}^{y(s)}\Big)\Big(\frac{\partial}{\partial s}  \Big)^{q+1}\Big((f^2(\xi)-L^2(s,\xi))^{q+\frac{1}{2}}\Big)d\xi\,
+V_1(s)=\Big(\frac{d}{ds}  \Big)^{q+1}\frac{const}{\sqrt{1+s^2}},
$$
where
$$
V_1(s)=\int\limits_{-x_o(1)}^{y_o(1)}\Big(\frac{\partial}{\partial s}  \Big)^{q+1}\Big((f_o^2(\xi)-L^2(s,\xi))^{q+\frac{1}{2}}\Big)d\xi.
$$
Note that, unlike it was for  the function $\Xi_1$ in the even-dimensional case, the function $V_1$ is well-defined only for $s\le 1$ and only if $ \|h\|_{C^1}$  is much less than $1$. Also,  even in that case, $V_1(s)$ is $C^{\infty}$ on $[0, 1)$ but not at $1$, where it is merely continuous.

Observe  that
$$
\Big(\frac{\partial}{\partial s}  \Big)^{q+1}\Big((f^2(\xi)-L^2(s,\xi))^{q+\frac{1}{2}}\Big)=\frac{J_1(s,\xi, f(\xi))}{\sqrt{f^2(\xi)-L^2(\xi)}},
$$
where $J_1(s,\xi,f)$ is some polynomial expression in $s$, $\xi$, $f$, $h(s)$, and the derivatives of $h$ at $s$.

Making the change of variables $\xi=-x(\sigma)$ in the  integral $\int_{-x(s)}^{-x_o(1)}$, and $\xi=y(\sigma)$  in the integral $\int_{y_o(1)}^{y(s)}$, we can rewrite the sum of integrals on the left hand side of (\ref{ass+}) as
$$
-\,\int\limits_s^{1}\left[\frac{J_1(s,-x(\sigma), L(\sigma,-x(\sigma)))}{\sqrt{L^2(\sigma, -x(\sigma))-L^2(s, -x(\sigma))}}\,\frac{dx}{ds}(\sigma)+
\frac{J_1(s,y(\sigma), L(\sigma, y(\sigma)))}{\sqrt{L^2(\sigma, y(\sigma))-L^2(s, y(\sigma))}}\,\frac{dy}{ds}(\sigma)\right]
d\sigma.
$$
Now write
$$
L^2(\sigma,\xi)-L^2(s, \xi)=(L(\sigma,\xi)-L(s, \xi))(L(\sigma,\xi)+L(s, \xi)),
$$
and
$$
L(\sigma,\xi)-L(s, \xi)=\sigma \xi+h(\sigma)-s\xi-h(s)=
(\sigma-s)(\xi+H(s,\sigma)),
$$
where
$$
H(s,\sigma)=\frac{h(\sigma)-h(s)}{\sigma-s}=\int\limits_0^1h'(s+(\sigma-s)\tau)d\tau
$$
is an infinitely smooth function of $s$ and $\sigma$.
Denote
$$
K_1(s,\sigma, \xi)=\frac{J_1(s, \xi, L(\sigma, \xi))}{\sqrt{ (\xi+ H(s,\sigma))(L(\sigma,\xi)+L(s, \xi))   }   }.
$$
The function $K_1$ is well-defined and infinitely smooth for all $s$, $\sigma$, $\xi$ satisfying
$ (\xi+ H(s,\sigma))$$(L(\sigma,\xi)+L(s, \xi))>0$. If $\|h\|_{C^1}$ is small enough, this condition is fulfilled whenever  $s$, $\sigma$$\in[\frac{1}{2},1]$ and $|\xi|>\frac{1}{2}$.

Now we can rewrite equation (\ref{ass+}) in the form
\begin{equation}\label{t1}
-\,\int\limits_s^{1}\Big(K_1(s,\sigma, -x(\sigma))
\,\frac{dx}{ds}(\sigma)\,+\,K_1(s,\sigma, y(\sigma))\,
\frac{dy}{ds}(\sigma)\Big)\,\frac{d\sigma}{\sqrt{\sigma-s}}=
\end{equation}
$$
-V_1(s)+\Big(\frac{d}{ds}  \Big)^{q+1}\frac{const}{\sqrt{1+s^2}}.
$$

Similarly, we can differentiate
 (\ref{dodd-}) and transform the resulting equation
into
\begin{equation}\label{t2}
-\,\int\limits_s^{1}\Big(K_2(s,\sigma, -x(\sigma))
\,\frac{dx}{ds}(\sigma)\,+\,K_2(s,\sigma, y(\sigma))\,
\frac{dy}{ds}(\sigma)\Big)\,\frac{d\sigma}{\sqrt{\sigma-s}}=-V_2(s),
\end{equation}
where $K_2$ is well-defined and infinitely smooth in the same range as $K_1$.
%
The function $V_2$ on the right hand side  of (\ref{t2})   is given by
$$
V_2(s)=\int\limits_{-x_o(1)}^{y_o(1)}\Big(\frac{\partial}{\partial s}  \Big)^{q}\Big((f_o^2(\xi)-L^2(s,\xi))^{q-\frac{1}{2}}L(s,\xi)\Big)d\xi,
$$
and everything that we said about $V_1$ applies to $V_2$ as well.

Equations (\ref{t1}) and (\ref{t2}) together can be written in the form
\begin{equation}\label{arik}
\int\limits_s^{1}\frac{K(s,\sigma,z(\sigma),\frac{dz}{ds}(\sigma))}{\sqrt{\sigma-s}}d\sigma=R(s),
\end{equation}
where, for $z=\left(
\begin{array}{cc}
x\\
y
\end{array}
\right),\,\,z'=\left(
\begin{array}{cc}
x'\,\\
y'\,
\end{array}
\right)\in\R^2$,
$$
K(s,\sigma,z,z')=\,-\,\left(
\begin{array}{cc}
K_1(s,\sigma, -x)
\,x'\,+\,K_1(s,\sigma, y)\,
y'\,\\
K_2(s,\sigma, -x)
\,x'\,+\,K_2(s,\sigma, y)\,
y'\,
\end{array}
\right),
$$
$$
R(s)=\left(
\begin{array}{cc}
-V_1(s)+\Big(\frac{d}{ds}  \Big)^{q+1}\frac{const}{\sqrt{1+s^2}}\,\\[10pt]
-V_2(s)\,
\end{array}
\right).
$$

By Lemma \ref{invert22} with $b=1$ (see Appendix),   equation (\ref{arik}) is equivalent to
\begin{equation}\label{dvi}
-G_2(s,s,z,z')+\int\limits_s^{1}\frac{\partial}{\partial s}G_2(s,\sigma, z(\sigma),\frac{dz}{ds}(\sigma))d\sigma=\widetilde{R}(s),
\end{equation}
where
$$
G_2(s,\sigma,z,z')=\int\limits_0^1\frac{K(s+\tau(\sigma-s),\sigma,z,z')}{\sqrt{\tau(1-\tau)}}d\tau,\qquad   \widetilde{R}(s)= \frac{d}{ds}\int\limits_s^{1}\frac{R(s')}{\sqrt{s'-s}}ds'.
$$
Note that
$$
G_2(s,s,z,z')=C\,\cdot\,K(s,s,z,z'),\qquad C=\int\limits_0^1\frac{d\tau}{\sqrt{\tau(1-\tau)}}.
$$
To reduce the resulting system of integro-differential equations to a pure system of integral equations we add two independent unknown functions $x'$, $y'$, denote
$z'=\left(
\begin{array}{cc}
x'\\
y'
\end{array}
\right)$, $z_o(s)=\left(
\begin{array}{cc}
x_o(s)\\
y_o(s)
\end{array}
\right)$, and
add
two new relations
$$
z(s)=-\int\limits_s^{1}z'(\sigma)d\sigma +z_o(1).
$$
Together with  (\ref{dvi}), they  lead to the system
\begin{equation}\label{hah44}
\mathbf{G}(s,Z(s))=\int\limits_s^{1}\mathbf{\Theta}(s,\sigma, Z(\sigma))d\sigma+\mathbf{\Xi}(s),\qquad Z=\left(
\begin{array}{cc}
z\\
z'
\end{array}
\right)=
\left(
\begin{array}{cccc}
x\\
y\\
x'\\
y'
\end{array}
\right).
\end{equation}
Here
$$
\mathbf{G}(s,Z)=\left(
\begin{array}{cc}
z\\
-G_2(s,s,z,z')\\
\end{array}
\right),\qquad
\mathbf{\Theta}(s,\sigma, Z)=-\,
\left(
\begin{array}{cc}
z'\\
\frac{\partial}{\partial s}G_2(s,\sigma,z,z')\\
\end{array}
\right),
$$
and
$$
\mathbf{\Xi}(s)=\left(
\begin{array}{cc}
z_o(1)\\
\widetilde{R}(s)
\end{array}
\right).
$$

In what follows, we will choose
 $h$  so that   $\|h\|_{C^1}$ is much less than $1$.
In this case,
${\mathbf G}$, ${\mathbf \Theta}$ are well-defined and infinitely smooth whenever $s$, $\sigma$ $\in [\frac{1}{2},1]$, $|x|$, $|y|$$>\frac{1}{2}$, $z'\in \R^2$, and $
{\mathbf \Xi}$ is well-defined  and infinitely smooth on $[\frac{1}{2},1)$.
Observe also  that
$$
D_Z{\mathbf G}\Big|_{(s,Z(s))}=
\left(
\begin{array}{cc}
\mathbf{I}&0\\
\mathbf{*}&\mathbf{A}\\
\end{array}
\right),
$$
where
$$
{\mathbf I}=\left(
\begin{array}{ccc}
1 & 0\\
0&
1
\end{array}
\right),\qquad
{\mathbf A(s,z)}=C\,\cdot\,{\mathbf E(s,z)},
$$
and
$$
{\mathbf E(s,z)}=
\left(
\begin{array}{cc}
K_1(s,s, -x)\,&\,
K_1(s,s, y)\,
\\
K_2(s,s, -x)
\,&\,\,K_2(s,s, y)\,
\end{array}
\right).
$$

The function
$$
Z_o(s)=\left(
\begin{array}{cc}
z_o(s)\\[3pt]
\frac{dz_o}{ds}(s)
\end{array}
\right)
=
\left(
\begin{array}{cccc}
x_o(s)\\
y_o(s)\\[3pt]
\frac{dx_o}{ds}(s)\\[3pt]
\frac{dy_o}{ds}(s)
\end{array}
\right)
$$
solves the system (\ref{hah44}) with ${\mathbf G}$, ${\mathbf \Theta}$, ${\mathbf \Xi}$ corresponding to $h\equiv 0$ (we will denote them by
${\mathbf G_o}$, ${\mathbf \Theta_o}$, ${\mathbf \Xi_o}$) on $[\frac{1}{2},1]$, say.

We claim that
\begin{equation}\label{detH111}
\det\,\Big(D_Z{\mathbf G_o}\Big|_{(s,Z_o(s))}  \Big)=\det({\mathbf A_o(s,z_o(s))})\neq 0\qquad \textrm{for all}\,\,\,s\in[\tfrac{1}{2},1].
\end{equation}
Indeed,
since $K_{1,2}(s,s,\xi)$  have the same signs as $J_{1,2}(s,\xi,L(s,\xi))$ and  since
$$
J_1(s,\xi, L(s,\xi))=(2q+1)!!\,\Big(-L(s,\xi)\frac{\partial}{\partial s}L(s,\xi)\Big)^{q+1},
$$
$$
J_2(s,\xi, L(s,\xi))=(2q-1)!!\,\Big(-L(s,\xi)\frac{\partial}{\partial s}L(s,\xi)\Big)^{q}L(s,\xi),
$$
we conclude that the matrix ${\mathbf A_o}(s,z_o(s))$ has the same sign pattern as the matrix
$$
\left(
\begin{array}{ccc}
(-1)^{q+1}& (-1)^{q+1}\\
(-1)^{q}(-x_o(s)) &  (-1)^{q}y_o(s)
\end{array}
\right),
$$
i.e.,  the signs in the first row are the same, and the signs in the second one are opposite.

Thus, (\ref{detH111}) follows.
In particular,
$$
\det\,\Big(D_Z{\mathbf G_o}\Big|_{(1,Z_o(1))}  \Big)\neq 0.
$$

Lemma \ref{dif} implies then that we can choose some small $\delta>0$ and construct a $C^k$-close to $Z_o(s)$ solution $Z(s)$ of (\ref{hah44}) on $[1-3\delta,1]$ whenever
${\mathbf G}$, ${\mathbf \Theta}$, ${\mathbf \Xi}$ are sufficiently close to ${\mathbf G_o}$, ${\mathbf \Theta_o}$, ${\mathbf \Xi_o}$ in $C^k$ on  certain compact sets. Since
$\mathbf{G}$, $\mathbf{\Theta}$, $\mathbf{\Xi}$ and their derivatives are   (integrals of) some explicit elementary expressions  in $Z$, $s$, $\sigma$, $h(s)$, and the derivatives of $h(s)$,
this closeness assumption will hold if $h$ is sufficiently close to zero with sufficiently many derivatives.
Moreover, since $h$ vanishes on $[1-\delta,1]$, the assumptions of Remark \ref{r} are satisfied and we have $Z(s)=Z_o(s)$ on $[1-\delta,1]$.

The $x$ and $y$ components of $Z$ solve the equations obtained by differentiating (\ref{dodd+}) and (\ref{dodd-}). The passage to  (\ref{dodd+}),  (\ref{dodd-}) is now exactly the same as in the even case.

\medskip

{\bf Step 3}.  From now on,  we change the point of view and switch to  the functions $R(\alpha)$ and $r(\alpha)$, $\alpha\in (0,\frac{\pi}{2})$, related
to $f$ by (\ref{f}). The functions $x$ and $y$, which  we have already constructed, implicitly define  $C^{\infty}$-functions $R_h(\alpha)$ and $r_h(\alpha)$ for all $\alpha$ with $\tan\alpha>1-3\delta$.

Instead of parameterizing hyperplanes by the slopes $s$ of the corresponding linear functions, we will
parameterize them by the angles $\beta$ they make with the $x_1$-axis, where $\beta$ is related to $s$ by $\tan\beta=s$.

Our next task will be to derive the equations that would ensure that all central sections corresponding to angles $\beta$ with $\tan\beta<1-2\delta$ are maximal
and of constant volume.
Note that those sections are already defined and satisfy these properties when
 $\tan \beta\in (1-3\delta,1-2\delta)$.

We remind that  if the volume of the central section  $K\cap v^{\perp}$ of a convex body $K$ is constant, then
\begin{equation}\label{mihalna1}
\frac{1}{d-1}({\mathcal R}\rho^{d-1}_K)(v)=\mathrm{vol}_{d-1}(K\cap v^{\perp})=const,
\end{equation}
where ${\mathcal R}$ is the {\it Spherical Radon Transform}, and $\rho_K$ is  the radial function of the body $K$.

Since
$$
\mathrm{vol}_{d-1}(K\cap (v^{\perp}+te_1))=\mathrm{vol}_{d-1}((K-te_1)\cap v^{\perp}),
$$
the
 central section corresponding to a unit vector $v$ with $\langle v,e_1\rangle\neq 0$
 is of maximal volume if and only if
\begin{equation}\label{mihalna2}
({\mathcal R}(\rho^{d-2}_K(\cdot)\left.\frac{\partial}{\partial t}\right|_{t=0}\rho_{K-te_1}(\cdot))(v)=0.
\end{equation}
Here $e_1$ is the ort along the $x_1$-axis.
We need these equations to hold for all unit vectors $v=(\pm\sin\alpha,\dots)\in {\mathbb S}^{d-1}$ corresponding to the angles  $\alpha$ with $\tan\alpha<1-2\delta$.

Note that when $K=K_f$ is the  body of revolution we are constructing,
these equations hold if $\tan\alpha\in (1-3\delta, 1-2\delta)$, and the left hand sides      of   (\ref{mihalna1}), (\ref{mihalna2})   are
already defined on the cap
$$
\{v\in{\mathbb S}^{d-1}:\,v=(\pm\sin\alpha,\dots), \alpha\in [0,\tfrac{\pi}{2}],\,\tan\alpha \ge 1-3\delta   \}
$$
and are smooth even rotation invariant functions there.   We will denote these functions  by $\varphi_h$ and
$\psi_h$ correspondingly.

Now we put $\varphi_h(v)=const$ and $\psi_h(v)=0$ for $v=(\pm\sin\alpha,\dots)$, $\tan\alpha\in [0,1-2\delta]$. This definition
agrees with the one we already have when $\tan\alpha\in [1-3\delta,1-2\delta]$, so $\varphi_h$ and $\psi_h$ are even rotation invariant infinitely smooth functions
on the entire sphere.

\medskip

Recall that the values of ${\mathcal R}g(v)$ for all $v=(\pm\sin\alpha,\dots)$ with $\tan\alpha>1-3\delta$
are completely defined by the values of
the even  function $g(u)$ for all  $u=(\pm\cos\alpha,\dots)$, $\tan\alpha>1-3\delta$
and {\it for bodies of revolution} (but not in general) the converse is also true (see the explicit inversion formula in  \cite{Ga}, page 433, formula (C.17)).

\medskip

Since  the equation  ${\mathcal R}g=\widetilde{g}$ with even $C^{\infty}$ right hand side $\widetilde{g}$ is equivalent to
$$
\frac{g(v)+g(-v)}{2}={\mathcal R^{-1}}\widetilde{g}(v),
$$
we can rewrite our equations      (\ref{mihalna1}) and (\ref{mihalna2})          as
$$
\rho^{d-1}_K(u)+\rho^{d-1}_K(-u)=2(d-1)({\mathcal R^{-1}}\varphi_h)(u)
$$
and
$$
\rho^{d-2}_K(u)\Big(\left.\tfrac{\partial}{\partial t}\right|_{t=0}\rho_{K-te_1}\Big)(u)+
\rho^{d-2}_K(-u)\Big(\left.\tfrac{\partial}{\partial t}\right|_{t=0}\rho_{K-te_1}\Big)(-u)=
2(d-1)({\mathcal R^{-1}}\psi_h)(u).
$$
The already constructed part of $\rho_K$ satisfies these equations for $u=(\pm\cos\alpha,\dots)$ with $\tan\alpha>1-3\delta$.

Since the Spherical Radon Transform commutes with rotations and our initial $\rho_K$ was rotation invariant,  the  even functions
$2(d-1){\mathcal R^{-1}}\varphi_h(u)$, $2(d-1){\mathcal R^{-1}}\psi_h(u)$
 are rotation invariant as well  and can be written as $\Phi_h(\alpha)$ and $\Psi_h(\alpha)$, where $u=(\pm\cos\alpha,\dots)$, $\alpha\in[0,\frac{\pi}{2}]$.
Note that the mappings $h\mapsto \Phi_h$, $h\mapsto\Psi_h$ are continuous from $C^{k+d}$ to $C^k$, say. Thus,  for all $h$ sufficiently close to zero in $C^{k+d}$,
$\Phi_h$ and $\Psi_h$ will be close to $\Phi_0\equiv2$ and
 $\Psi_0\equiv 0$  in $C^k$.

We will be looking for a rotation invariant solution $\rho_K$,  which will be described in terms of two functions $R(\alpha)$ and $r(\alpha)$ related to it by (\ref{tryuk}).
Equation   (\ref{mihalna1})        translates into
\begin{equation}\label{+r}
r^{d-1}(\alpha)+R^{d-1}(\alpha)=\Phi_h(\alpha).
\end{equation}
To rewrite  equation       (\ref{mihalna2}) , observe that
\begin{equation}\label{ap}
\rho^{d-2}_K(u)\Big(\left.\tfrac{\partial}{\partial t}\right|_{t=0}\rho_{K-te_1}\Big)(u)+
\rho^{d-2}_K(-u)\Big(\left.\tfrac{\partial}{\partial t}\right|_{t=0}\rho_{K-te_1}\Big)(-u)=
\end{equation}
$$
-\Big[R^{d-3}(\alpha)(R(\alpha)\sin\alpha)'-r^{d-3}(\alpha)(r(\alpha)\sin\alpha)'\Big]
$$
(see Lemma \ref{trum} in Appendix). Thus,  equation  (\ref{mihalna2})  translates into
$$
R^{d-3}(\alpha)(R(\alpha)\sin\alpha)'-r^{d-3}(\alpha)(r(\alpha)\sin\alpha)'=
-\Psi_h(\alpha).
$$
Multiplying it by $(d-2)\sin^{d-3}\alpha$,  we obtain
\begin{equation}\label{+r1}
\Big((R(\alpha)\sin\alpha)^{d-2}-(r(\alpha)\sin\alpha)^{d-2}\Big)'=-(d-2)\Psi_h(\alpha)\sin^{d-3}\alpha.
\end{equation}
Taking into account the condition $R(\frac{\pi}{2})=r(\frac{\pi}{2})$ and integrating, we see that (\ref{+r1}) can also be written as
\begin{equation}\label{-r}
R^{d-2}(\alpha)-r^{d-2}(\alpha)=\frac{\Theta_h(\alpha)}{\sin^{d-2}\alpha},
\end{equation}
where
$$
\Theta_h(\alpha)=(d-2)\int\limits_{\alpha}^{\pi/2}\Psi_h(\beta)\sin^{d-3}\beta d\beta.
$$

Note that equations (\ref{+r}), (\ref{-r}) together with conditions $R(\alpha)>0$ and $r(\alpha)>0$ determine $R(\alpha)$ and $r(\alpha)$ uniquely, and the originally constructed functions $R_h$ and $r_h$ satisfy these equations for  all $\alpha\in[0,\frac{\pi}{2}]$    with $\tan\alpha\ge 1-3\delta$. Thus, any solution $R$,  $r$  of this system will satisfy $R(\alpha)=R_h(\alpha)$,
$r(\alpha)=r_h(\alpha)$ in this range.

It remains to show that the solutions $R$ and $r$ exist at all and  define a {\it convex body}  if $h$ is chosen appropriately. Note that when $h$ is small with several derivatives the functions $\Phi_h-2$ and $\Psi_h$ are close to zero uniformly with several derivatives. The only problem is that the right hand side of (\ref{-r})
can blow up as $\alpha\to 0^+$. To prevent it, we will choose the perturbation function $h$ close to $0$ in $C^{2d+k}$ so that
\begin{equation}\label{bu}
\Theta_h(0)=\Theta_h'(0)=...=\Theta_h^{(d+k-1)}(0)=0.
\end{equation}
Then, the right hand side of (\ref{-r}) will be close to zero in $C^{k}([0,\frac{\pi}{2}])$.
Since  the map
$$
\mathbf{D}:\,(R,r)\,\,\mapsto\,\,(
R^{d-1}+r^{d-1},\,
R^{d-2}-r^{d-2})
$$
is smoothly invertible near the point $(1,1)$  by the inverse function theorem,  the functions $R$, $r$ exist in this case on the entire interval $[0,\frac{\pi}{2}]$, and are close to $1$ in $C^2$. Moreover, $R'(0)=r'(0)=0$, because $\Phi^{'}_h(0)=0$ (otherwise the function ${\mathcal R^{-1}}\varphi_h$ would not be smooth at $(1,0,\dots,0)$), and
the right hand side of (\ref{-r}) is $o(\alpha)$ as $\alpha\to 0^+$. This is enough to ensure that the body given by $R$ and $r$ is convex and corresponds to some strictly
 concave function $f$ defined on $[-r(0),R(0)]$.

Finally, to prove the existence of a perturbation function $h$ for which cancellation conditions (\ref{bu}) hold, we use the Borsuk-Ulam Theorem, stating that  a continuous map from
${\mathbb S}^m$ to ${\mathbb R}^m$, taking the antipodal points to antipodal points, contains zero in its image.
For $x=(x_1,...,x_{d+k+1})\in {\mathbb S}^{d+k}$ we define $h_x=\sum\limits_{j=1}^{d+k+1}x_jh_j$, where $h_j$ are  not identically zero smooth
 functions on $[1-2\delta, 1-\delta]$ with pairwise disjoint supports.
We define the map $\mathbf{B}:\,{\mathbb S}^{d+k}\to {\mathbb R}^{d+k}$ by
$$
\mathbf{B}:\,(x_1,...,x_{d+k+1})\,\,\mapsto \,\,h_x\,\,\mapsto\,\,(\Theta_{h_x}(0),\Theta_{h_x}'(0),...,\Theta_{h_x}^{(d+k-1)}(0)).
$$
Observing  that  $R_{(-h)}=r_{h}$ and $r_{(-h)}=R_{h}$, and   using   the linearity of the Inverse Spherical Radon Transform, we conclude that
 the map
$h\mapsto \Psi_h$ is odd. Hence, $\mathbf{B}$  maps antipodal points to antipodal points, and  there exists some not identically zero $h$ for which   (\ref{bu}) holds.
This completes the proof of Theorem \ref{klee} in the case of an  odd dimension.

\section{Appendix }

All results collected in this appendix are well-known. However, since we targeted this article  not exclusively  at specialists in integral equations and since in many cases it was much harder to find an exact reference than to write a full proof, we decided to present them here.
The reader should also keep in mind that we tailored the exact statements to our particular needs, so we cut corners whenever possible to reduce the presentation to the bare minimum.

\bl\label{ll0}
Let  $g$, $q\in C([a,b])$ be two functions with values in $\R^m$. Then, $g\equiv q$ if and only if for every $s\in[a,b)$ the equality
$$
\int\limits_s^b\frac{g(s')}{\sqrt{s'-s}}\,ds'\,=\int\limits_s^b\frac{q(s')}{\sqrt{s'-s}}\,ds'
$$
holds.
\el
\bp
The only non-trivial statement here is that the equality of integrals implies the equality of functions. To prove it,  take $t<b$, multiply both parts by $\frac{1}{\sqrt{s-t}}$,
and integrate with respect to $s$ from $t$ to $b$.  Making the change of variables 
\linebreak
$s=t+\tau(s'-t)$, performing the integration with respect to $\tau$ first,  and canceling
the common factor
$\int_0^1\frac{d\tau}{\sqrt{\tau(1-\tau)}}>0$,
we obtain
$$
\int\limits_t^bg(s')\,ds'\,=\int\limits_t^bq(s')\,ds' \qquad \text{for all}\qquad t\in [a,b).
$$
The result follows now from the fundamental theorem of calculus and the continuity assumption.
\ep

\bl\label{ll1}
Let $R:[a,b]\to \R^m$. If $R\in C([a,b])\cap C^{\infty}([a,b))$, then the function
$$
s\,\mapsto\,\int\limits_s^b\frac{R(s')}{\sqrt{s'-s}}\,ds'
$$
belongs to $ C^{\infty}([a,b))$ and tends to zero as $s\to b^-$.
\el
\bp
The second statement follows from the crude bound
$$
\int\limits_s^b\frac{|R(s')|}{\sqrt{s'-s}}\,ds'\,\le\,2\|R\|_{C([a,b])}\sqrt{b-s}\,\to \,0
$$
as $s\to b^-$. To prove the first one,  fix $\delta>0$. For  $s<b-\delta$,  we have
$$
\int\limits_s^b\frac{R(s')}{\sqrt{s'-s}}\,ds'=\int\limits_s^{s+\delta}\frac{R(s')}{\sqrt{s'-s}}\,ds'\,+\,\int\limits_{s+\delta}^b\frac{R(s')}{\sqrt{s'-s}}\,ds'.
$$
The kernel
$\frac{1}{\sqrt{s'-s}}$ is $C^{\infty}$-smooth in $s$ for $s'>s+\delta$, so the second integral is  infinitely smooth on $[a,b)$ because $R$ is infinitely smooth there. The first integral can be rewritten as
$$
\int\limits_0^{\delta}\frac{R(s+\tau)}{\sqrt{\tau}}\,d\tau.
$$
Since $s+\tau$ stays away from $b$ when $s$ stays away from $b-\delta$, this integral  also defines a $C^{\infty}$ function on $[a,b-\delta)$. Since $\delta$ is arbitrary, the lemma follows.
\ep
\bl\label{ll111}
Let $U\in C([a,b]^2)$. Then the function
$$
s\,\mapsto\,\int\limits_s^b\frac{U(s,\sigma)}{\sqrt{\sigma-s}}\,d\sigma
$$
extended by zero to $s=b$ is continuous.
\el
\bp
Let $s<b$. Fix $\delta\in (0,b-s)$ and let $|s-s'|<\frac{\delta}{2}$.  Write
$$
\Big|\int\limits_{s^{'}}^b\frac{U(s',\sigma)}{\sqrt{\sigma-s'}}\,d\sigma \,-\,\int\limits_s^b\frac{U(s,\sigma)}{\sqrt{\sigma-s}}\,d\sigma\Big|\le
$$
$$
\int\limits_{s^{'}}^{s+\delta}\frac{|U(s',\sigma)|}{\sqrt{\sigma-s'}}\,d\sigma \,+\,\int\limits_s^{s+\delta}\frac{|U(s,\sigma)|}{\sqrt{\sigma-s}}\,d\sigma\,+\,\int\limits_{s+\delta}^{b}\Big|\frac{U(s',\sigma)}{\sqrt{\sigma-s'}}\,-\,\frac{U(s,\sigma)}{\sqrt{\sigma-s}}\Big|\,d\sigma.
$$
The first two integrals are bounded by $2\|U\|_C\sqrt{\delta+|s'-s|}\le 2\|U\|_C\sqrt{2\delta}$ and $2\|U\|_C\sqrt{\delta}$ respectively.
The third one tends to zero as $s'\to s$ because the integrand tends to zero uniformly. The continuity at $b$ follows  from the inequality
$$
\Big|\int\limits_s^b\frac{U(s,\sigma)}{\sqrt{\sigma-s}}\,d\sigma\Big|\le2\|U\|_C\sqrt{b-s}.
$$
\ep
\bl\label{invert22}
Let $U\in C^{\infty}([a,b]^2)$, and let
$R\in C([a,b])\cap C^{\infty}([a,b))$.  Then  the equation
\begin{equation}\label{hrum}
\int\limits_s^b\frac{U(s, \sigma)}{\sqrt{\sigma-s}}\,d\sigma\,=R(s)
\end{equation}
holds for all  $s\in[a,b)$
 if and only if so does  the equation
 $$
 -\,V(s,s)+\int\limits_s^b\frac{\partial V}{\partial s}(s,\sigma)d\sigma=\widetilde{R}(s),
 $$
 where
 $$
 V(s,\sigma)=\int\limits_0^1\frac{U(s+\tau(\sigma-s),\sigma)}{\sqrt{\tau(1-\tau)}}\,d\tau,\qquad \widetilde{R}(s)=\frac{d}{ds}\int\limits_s^b\frac{R(s')}{\sqrt{s'-s}}\,ds'.
 $$
\el
\bp
Observe that our assumption on $U$ implies that $V\in C^{\infty}([a,b]^2)$. Also,
by the previous lemma, $\widetilde{R}(s)$ is well-defined for $s\in [a,b)$. By Lemmata  \ref{ll111} and \ref{ll1} both parts of (\ref{hrum}) are continuous functions.  Therefore, by Lemma \ref{ll0}, equation
(\ref{hrum}) is equivalent to
$$
\iint\limits_{s<s'<\sigma<b}\frac{U(s', \sigma)}{\sqrt{\sigma-s'}\sqrt{s'-s}}\,ds' d\sigma=\int\limits_s^b\frac{R(s')}{\sqrt{s'-s}}\,ds'.
$$
Making the change of variables  $s'=t+\tau(\sigma-s)$ and performing the integration with respect to $\tau$ first,
we can rewrite the left hand side as
$\int_s^bV(s,\sigma)d\sigma$.
Observe that both parts tend to zero as $s\to b^-$, and are differentiable functions on $[a,b)$. Thus, their equality is equivalent to the equality of their derivatives.
\ep
\bl
{\rm(}\textbf{Banach Fixed Point Theorem}, \cite{Ba}{\rm)}. Suppose that $(X,d)$ is a complete metric space and $T$ is a mapping from $X$ to $X$ satisfying 
$$d(T(x),T(y))\le k\, d(x,y)$$ for
all $x,y\in X$ with
some $k\in(0,1)$.
Then

1) there exists a unique fixed point $z_*$ of the mapping $T$,

2) the sequence of Picard iterations $z_{k+1}=T(z_k)$ starting with any point $z_0\in X$ converges to $z_*$,

3) for every point $z\in X$, we have $d(z,z_*) \le\dfrac{d(T(z),z)}{1-k}$.
\el

\bl\label{dif}
Let $\Omega$ be a domain in $\R^m$,  let  $[a,b]\subset \R$, and let  ${\mathbf G_o}:\,[a,b]\times \Omega\to \R^m$, ${\mathbf \Theta_o}:\,[a,b]^2\times \Omega\,\to \R^m$,  ${\mathbf \Xi_o}:\, [a,b)\to \R^m$,  $Z_o:\, [a,b]\to \Omega$. Assume  that ${\mathbf G_o}$, ${\mathbf \Theta_o}$, $Z_o$  are infinitely smooth,  ${\mathbf \Xi_o}$ is continuous, and
$$
{\mathbf G_o}(s,Z_o(s))=\int\limits_s^{b}{\mathbf \Theta_o}(s,\sigma, Z_o(\sigma))d\sigma+{\mathbf \Xi_o}(s)
$$
for all $s\in [a,b]$.

Then ${\mathbf \Xi_o}$ extends to a $C^{\infty}$-function on $[a,b]$.
If, in addition,
$$
 \det\Big(D_Z{\mathbf G_o}\Big|_{(b,Z_o(b))}   \Big)\neq 0,
$$
then there exist
 $ \varepsilon>0$, $\delta>0$, such that    on the interval $[b-3\delta,b]$,  every perturbed equation
\begin{equation}\label{purd}
{\mathbf G}(s,Z(s))=\int\limits_s^{b}{\mathbf \Theta}(s,\sigma, Z(\sigma))d\sigma+{\mathbf \Xi}(s)
\end{equation}
 has a unique continuous solution $Z(s)$  satisfying $\| Z-Z_o(b)\|_{C([b-3\delta,b])}<\varepsilon$, provided
 ${\mathbf G}$, ${\mathbf \Theta}$, ${\mathbf \Xi}$ are infinitely smooth   and
$$
\|{\mathbf G}-{\mathbf G_o}\|_{C^1([b-3\delta,b]\times B)},\qquad \|{\mathbf \Theta}-{\mathbf \Theta_o}\|_{C^1([b-3\delta,b]^2\times B)},\qquad
\|{\mathbf \Xi}-{\mathbf \Xi_o}\|_{C^1([b-3\delta,b])}
$$
are small enough. This solution is infinitely smooth and close to $Z_o$ in $C^k$, provided that
$$
\|{\mathbf G}-{\mathbf G_o}\|_{C^k([b-3\delta,b]\times B)},\qquad \|{\mathbf \Theta}-{\mathbf \Theta_o}\|_{C^k([b-3\delta,b]^2\times B)},\qquad
\|{\mathbf \Xi}-{\mathbf \Xi_o}\|_{C^k([b-3\delta,b])}
$$
are small enough. Moreover, the solutions corresponding to two different triples  ${\mathbf G}$, ${\mathbf \Theta}$, ${\mathbf \Xi}$ that are close in corresponding
$C^k$ are $C^k$ close to each other.

Here $B$ is the closed ball of radius $\varepsilon$ centered at $Z_o(b)$.
\el
\bp
The first statement is obvious because  all terms in the unperturbed equation except ${\mathbf \Xi_o}$ are defined and infinitely smooth on the entire interval $[a,b]$.

Next, denote ${\mathbf Q}=D_Z{\mathbf G_o}(b,Z_o(b))$ and observe that there exist $\varepsilon>0$, $\delta_1>0$ such that   $B\subset \Omega$,
 and for all $s\in [b-3\delta_1,b]$  and all $Z$ such that $|Z-Z_o(b)|<\varepsilon$,  we have
$$
|Z_o(s)-Z_o(b)|<\frac{\varepsilon}{8}\quad\textrm{and}\quad\|D_Z{\mathbf G_o}((s,Z))-{\mathbf Q}\|\le \frac{1}{8\|{\mathbf Q^{-1}}\|}.
$$
The perturbed equation (\ref{purd}) can be rewritten as
$Z(s)=(TZ)(s)$,
where
$$
(TZ)(s)=Z(s)-{\mathbf Q^{-1}}\Big[{\mathbf G}(s,Z(s))-\int\limits_s^{b}{\mathbf \Theta}(s,\sigma, Z(\sigma))d\sigma-{\mathbf \Xi}(s)\Big].
$$
We will show that if $\delta\in(0,\delta_1)$ is small enough and $\|{\mathbf G}-{\mathbf G_o}\|_{C^1([b-3\delta,b])}\le \dfrac{1}{8\|{\mathbf Q^{-1}}\|}$,
then  $T$  is a contraction on the set $X$ of continuous functions mapping $C([b-3\delta,b])$ to
$B$.  Take two functions $Z_1$ and $Z_2$ in $C([b-3\delta,b])$ with values in $B$ and notice that
$$
\Big|\int\limits_s^{b}\Big({\mathbf \Theta}(s,\sigma, Z_1(\sigma))-{\mathbf \Theta}(s,\sigma, Z_2(\sigma))\Big)d\sigma\Big|\le
$$
$$
 3\delta\Big[\max\limits_{s',\sigma'\in[b-3\delta,b], |Z-Z_o(b)|\le \varepsilon}\|D_Z{\mathbf \Theta}(s',\sigma', Z) \|\Big]\|Z_1-Z_2\|_{C([b-3\delta,b])}
$$
for all $s\in [b-3\delta,b]$.  So, if
$\delta$ is chosen so small that 
$$
3\delta\,\|{\mathbf Q^{-1}}\|\, \left(\|{\mathbf \Theta_o} \|_{C^1([b-3\delta,b])}+1      \right)\le\frac{1}{4}
$$
and if $\|{\mathbf \Theta}-{\mathbf \Theta_o} \|_{C^1([b-3\delta,b])}<1$, we have
$$
3\delta\,\|{\mathbf Q^{-1}}\|\,\max\limits_{s',\sigma'\in[b-3\delta,b], |Z-Z_o(b)|\le \varepsilon}\|D_Z{\mathbf \Theta} (s',\sigma', Z)\|\le
$$
$$
 3\delta\, \|{\mathbf Q^{-1}}\|\Big(\|{\mathbf \Theta_o} \|_{C^1([b-3\delta,b])}+
\|{\mathbf \Theta}-{\mathbf \Theta_o} \|_{C^1([b-3\delta,b])}  \Big)<\frac{1}{4},
$$
and
$$
\Big|{\mathbf Q^{-1}}\Big(\int\limits_s^{b}{\mathbf \Theta}(s,\sigma, Z_1(\sigma))d\sigma- \int\limits_s^{b}{\mathbf \Theta}(s,\sigma, Z_2(\sigma))d\sigma  \Big)\Big|\le \frac{1}{4}\|Z_1-Z_2\|_{C([b-3\delta,b])}
$$
for all $s\in [b-3\delta,b]$.

Consider the function ${\mathbf H}(s, Z)=Z-{\mathbf Q^{-1}}{\mathbf G}(s,Z)$. Note that
$$
\| D_Z{\mathbf H} \|=\|{\mathbf I}-{\mathbf Q^{-1}}D_Z{\mathbf G}\|=\|{\mathbf Q^{-1}}(D_Z{\mathbf G}-{\mathbf Q})\|\le
$$
$$
 \|{\mathbf Q^{-1}}   \| \Big[ \| D_Z{\mathbf G_o}-{\mathbf Q}  \|+\| D_Z{\mathbf G}-{\mathbf G_o}  \|\Big]\le\frac {1}{4}
$$
when $s\in [b-3\delta,b]$ and $|Z-Z_o(b)|\le \varepsilon$.
Thus,
$$
\Big|\Big(Z_1(s)-{\mathbf Q^{-1}}{\mathbf G}(s,Z_1(s))\Big) -  \Big(Z_2(s)-{\mathbf Q^{-1}}{\mathbf G}(s,Z_2(s))\Big) \Big|=\Big|{\mathbf H}(s, Z_1)-{\mathbf H}(s, Z_2)   \Big|\le
$$
$$
 \Big[\max\limits_{s'\in [b-3\delta,b],\,|Z-Z_o(b)|\le \varepsilon }\| D_Z {\mathbf H}(s',Z)\|\Big] \|Z_1-Z_2\|_{C([b-3\delta,b])} \le \frac{1}{4}\|Z_1-Z_2\|_{C([b-3\delta,b])}.
$$
Bringing these estimates together, we see that
$$
\Big|(TZ_1)(s)-(TZ_2)(s)\Big|\le\frac{1}{2}\|Z_1-Z_2\|_{C([b-3\delta,b])}
$$
for all $s\in [b-3\delta,b]$.  To apply the Banach fixed point theorem it remains to show that $T$ maps
$X$ to itself. To this end,
notice that
$$
\Big|(TZ_o)(s)-Z_o(s)   \Big|=\Big| {\mathbf Q^{-1}} \Big[\Big({\mathbf G}(s,Z_o(s)) -  {\mathbf G_o}(s,Z_o(s))  \Big)  -
$$
$$
\int\limits_s^{b}\Big( {\mathbf  \Theta}(s,\sigma, Z_o(\sigma))-{\mathbf  \Theta_o}(s,\sigma, Z_o(\sigma))\Big)d\sigma-\Big({\mathbf \Xi(s)}-{\mathbf \Xi_o}(s)\Big)\Big] \,
\Big|\le
$$
\begin{equation}\label{hryu11}
\| {\mathbf Q^{-1}} \| \, \Big( \|{\mathbf G} -  {\mathbf G_o} \|_{C([b-3\delta,b]\times B)}  +
\end{equation}
$$
3\delta\|{\mathbf \Theta}-{\mathbf \Theta_o}\|_{C([b-3\delta,b]^2\times B)} + \|{\mathbf \Xi(s)}-{\mathbf \Xi_o}(s)   \|_{C([b-3\delta,b])}\Big)<\frac{\varepsilon}{4},
$$
provided that the $C$-norms in the last two lines are small enough.

Let $Z\in X$. It is obvious that $TZ$ is a continuous function on $[b-3\delta,b]$. Also,
$$
|Z(s)-Z_o(s)|\le|Z(s)-Z_o(b)|+|Z_o(s)-Z_o(b)|\le \varepsilon+\frac{\varepsilon}{8}=\frac{5\varepsilon}{8}
$$
for all $s\in [b-3\delta,b]$, so $\|Z-Z_o\|_{C([b-3\delta,b])}\le\frac{5\varepsilon}{8}$, and
$$
|TZ(s)-Z_o(b)|\le|(TZ)(s)-(TZ_o)(s)|+|(TZ_o)(s)-Z_o(s)|+ |Z_o(s)-Z(b)|   \le
$$
$$
\frac{1}{2}\cdot\frac{5\varepsilon}{8} +\frac{\varepsilon}{8}+\frac{\varepsilon}{8}<\varepsilon.
$$
Thus, $TZ\in X$ as well. This completes the proof of the existence and uniqueness part of the lemma.

To show that $Z$ is smooth, notice that the right hand side of equation (\ref{purd}) is a $C^1$ function for every $Z\in X$.
Thus, the left hand side  ${\mathbf G}(s,Z(s))$ is also $C^1$. Since  ${\mathbf G}\in C^{\infty}([a,b]\times\Omega)$ and $D_Z{\mathbf G}(s,Z)\neq 0$ whenever $s\in [b-3\delta, b]$, $|Z-Z_o(b)|<\varepsilon$, we conclude by the implicit function theorem that  $Z\in C^1$ and, moreover,
$$
\frac{dZ}{ds}(s)=
$$
$$
\Big(D_Z{\mathbf G}(s,Z(s))\Big)^{-1}\Big(-\frac{\partial {\mathbf G}}{\partial s}(s,Z(s))- {\mathbf \Theta}(s,s, Z(s)+\int\limits_s^{b}\frac{\partial {\mathbf \Theta}}{\partial s}(s,\sigma, Z(\sigma))d\sigma+\frac{d{\mathbf \Xi}}{ds}(s)  \Big).
$$
Differentiating this identity again and again and plugging the expression for the derivative $\frac{dZ}{ds}$ into the right hand side after every differentiation, we see that $Z$ is infinitely smooth, and, moreover, $(\frac{d}{ds})^kZ$ can be written as some explicit expression involving only $Z$ itself and various partial derivatives
of the functions ${\mathbf G}$, ${\mathbf \Theta}$, ${\mathbf \Xi}$ of orders up to $k$. We see from here that to show that $Z$ is close to $Z_o$ in $C^k$ under the condition
that the norms
$$
\|{\mathbf G}-{\mathbf G_o}\|_{C^k([b-3\delta,b]\times B)},\qquad \|{\mathbf \Theta}-{\mathbf \Theta_o}\|_{C^k([b-3\delta,b]^2\times B)},\qquad
\|{\mathbf \Xi}-{\mathbf \Xi_o}\|_{C^k([b-3\delta,b])}
$$
are small, it suffices to show that, under this condition, the norm $\|Z-Z_o\|_{C([b-3\delta,b])}$ is small. By the third part of the Banach fixed point theorem, this would
follow from the smallness of $\|TZ_o-Z_o\|_{C([b-3\delta,b])}$. But we have already estimated this difference  by
$$
\| {\mathbf Q^{-1}} \| \, \Big( \|{\mathbf G} -  {\mathbf G_o} \|_{C([b-3\delta,b]\times B)}  +3\delta\|\Theta-\Theta_o\|_{C([b-3\delta,b]^2\times B)} + \|\Xi(s)-\Xi_o(s)   \|_{C([b-3\delta,b])}\Big)
$$
in (\ref{hryu11}).

Exactly the same argument can be used to prove the last statement of the lemma.
\ep
\begin{remark}\label{r}
If ${\mathbf \Xi}={\mathbf \Xi_o}$ on $[b-\delta, b]$, then to check that $\|{\mathbf \Xi}-{\mathbf \Xi_o}\|_{C^k([b-3\delta,b])}$ is small, it suffices to check that
$\|{\mathbf \Xi}-{\mathbf \Xi_o}\|_{C^k([b-3\delta,b-\delta])}$ is small. If, in addition, ${\mathbf  G}(s,Z)={\mathbf  G_o}(s,Z)$ for all $s\in[b-\delta,b]$, and ${\mathbf  \Theta}(s,\sigma, Z)={\mathbf  \Theta_o}(s,\sigma, Z)$ for all $s,\sigma\in [b-\delta,b]$, then the solution $Z$, whose existence and uniqueness is asserted in Lemma \ref{dif},
 coincides with $Z_o$ on $[b-\delta,b]$.
\end{remark}
This follows  from the fact that if $Z=Z_o$ on $[b-\delta,b]$, then $TZ = Z$ on $[b-\delta,b]$ as well, so
if
we start the Picard iterations with $Z_o$,   the values on this interval will never  change.

\bl\label{trum}
Let $K$ be a body of revolution around the $x_1$-axis and let $\rho_K$ be the radial function of $K$. Then
$$
\rho^{d-2}_K(u)\Big(\left.\tfrac{\partial}{\partial t}\right|_{t=0}\rho_{K-te_1}\Big)(u)+
\rho^{d-2}_K(-u)\Big(\left.\tfrac{\partial}{\partial t}\right|_{t=0}\rho_{K-te_1}\Big)(-u)=
$$
$$
-\Big[R^{d-3}(\alpha)(R(\alpha)\sin\alpha)'-r^{d-3}(\alpha)(r(\alpha)\sin\alpha)'\Big]
$$
with $R$ and $r$  defined by
$$
\rho_K(u)= \begin{cases}
R(\alpha)\qquad &\text{if \,}u_1>0,\\
r(\alpha) & \text{if \,}u_1<0,
\end{cases}
$$
where  $u=(u_1,\dots)\in {\mathbb S}^{d-1}$ and $\alpha\in(0,\frac{\pi}{2})$, $\cos\alpha=|u_1|$.
\el
\bp
Denote by $W$ the $(x_1,x_2)$-plane. Let
 $l$ be the line $\{(x_1,x_2,\dots)\in W:\,x_2=x_1\tan \alpha \}$, where $\alpha\in (0,\pi/2)$. For a small $t>0$ we denote  by $l_{t}$  the line $\{x\in W:\,x_2=(x_1-t)\tan\alpha\}$.

\pagebreak

\begin{figure}[ht]
\includegraphics[width=250pt]{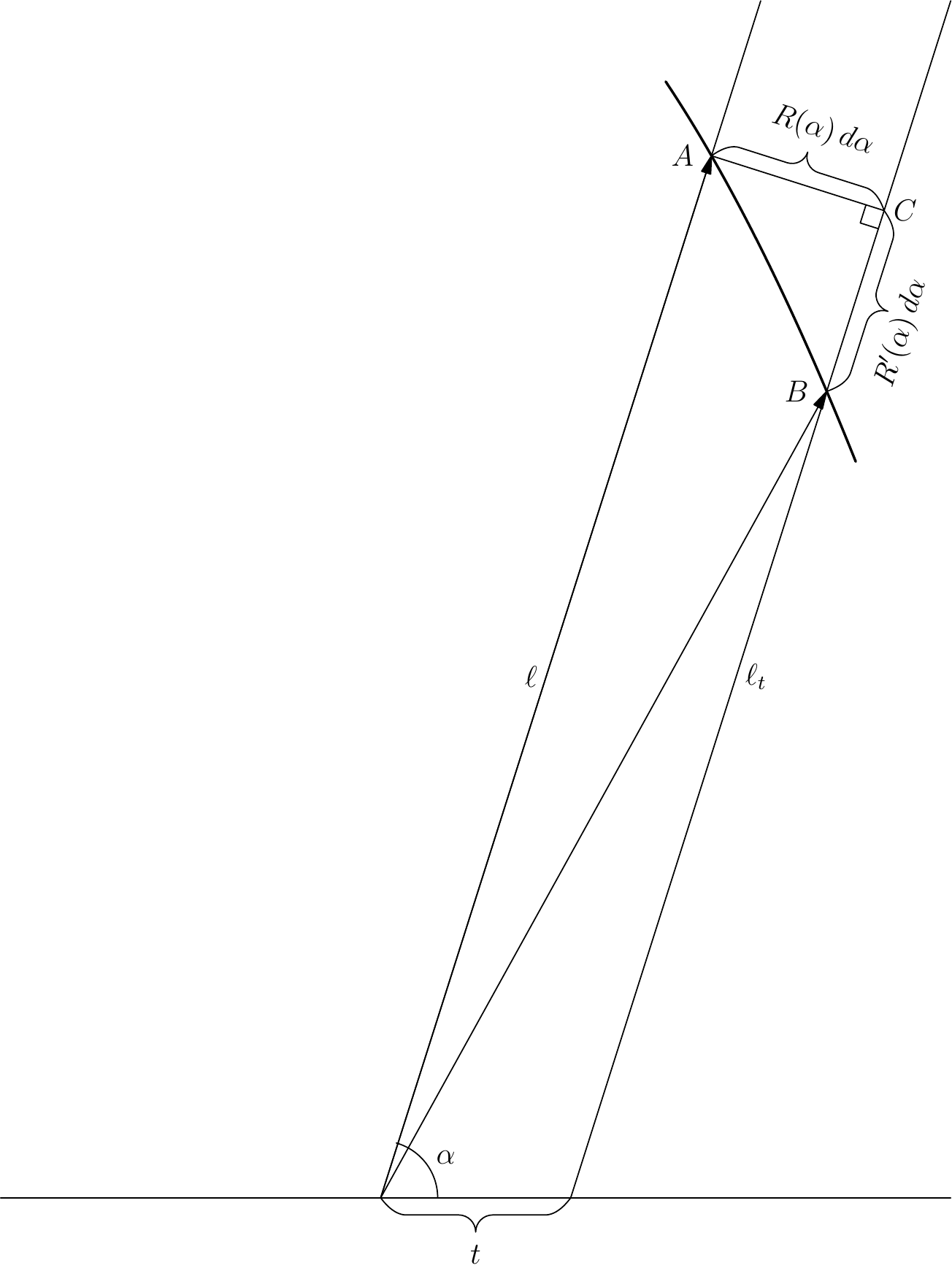}
\caption{$\rho_{K-te_1}(u)=R(\alpha)-t\cos\alpha-t\frac{R^{'}(\alpha)}{R(\alpha)}\sin\alpha$}
\label{pic6}
\end{figure}

Denote by $A$ and   $B$
 the ``top" points of intersection of  the boundary of $K$
with  $l$ and $l_{t}$  correspondingly.  Let $C$ be the point of intersection of $l_{t}$ with the hyperplane   orthogonal  to $l$  and passing through $A$ (see Figure \ref{pic6}).
Observe  that $A$, $B$, $C\in W$ and that  $K\cap l_{t}$ is the one-dimensional {\it central section} of the shifted body $K-te_1$.

By elementary geometry,
$$
\rho_{K-te_1}(u)=R(\alpha)-t\cos\alpha-t \sin \alpha\,\tan CAB,
$$
where $u=(\cos\alpha,\dots)\in {\mathbb S}^{d-1}$ and $\alpha\in(0,\frac{\pi}{2})$.
Observe that, up to  terms of order $t^2$, we have
$$
\tan CAB=\frac{R'(\alpha)}{R(\alpha)}.
$$
Hence,
$$
\rho_{K-te_1}(u)=R(\alpha)-t\cos\alpha-t \frac{R'(\alpha)}{R(\alpha)}\sin \alpha+o(t^2).
$$
Similarly,
$$
\rho_{K-te_1}(-u)=r(\alpha)+t\cos\alpha+t \frac{r'(\alpha)}{r(\alpha)}\sin \alpha +o(t^2).
$$
Finally,
$$
\rho^{d-2}_K(u)\Big(\left.\tfrac{\partial}{\partial t}\right|_{t=0}\rho_{K-te_1}\Big)(u)\,+\,
\rho^{d-2}_K(-u)\Big(\left.\tfrac{\partial}{\partial t}\right|_{t=0}\rho_{K-te_1}\Big)(-u)=
$$
$$
R^{d-2}(\alpha)(-\cos\alpha-\frac{R'(\alpha)}{R(\alpha)}\sin\alpha)-
r^{d-2}(\alpha)(-\cos\alpha- \frac{r'(\alpha)}{r(\alpha)}\sin \alpha)=
$$
$$
R^{d-3}(\alpha)(-R(\alpha)\cos\alpha-R'(\alpha)\sin\alpha)-
r^{d-3}(\alpha)(-r(\alpha)\cos\alpha- r'(\alpha)\sin \alpha)=
$$
$$
-\Big[R^{d-3}(\alpha)(R(\alpha)\sin\alpha)'-r^{d-3}(\alpha)(r(\alpha)\sin\alpha)'\Big].
$$
\ep

\end{document}